\documentclass[a4paper,11pt]{article}
\usepackage{authblk}
\usepackage[utf8]{inputenc}
\usepackage{amsmath,amsthm}
\usepackage{soul}
\usepackage[normalem]{ulem}
\usepackage{amssymb,latexsym}

\usepackage{comment}
\usepackage{caption}
\usepackage{mathrsfs}

\numberwithin{equation}{section}
\usepackage[
left=34mm,
right=34mm,
top=25mm,
bottom=25mm,
includeheadfoot,
]{geometry}

\usepackage{xcolor}

\theoremstyle{plain}
\newtheorem{thm}{Theorem}[section]
\newtheorem{prop}{Proposition}[section]
\newtheorem{defi}{Definition}[section]

\newtheorem{remq}{Remark}[section]
\newtheorem{lem}{Lemma}[section]

\numberwithin{equation}{section} 
\newcommand{\R}{\mathbb{R}}

\newcommand{\iy}{\infty}
\newcommand{\p}{\partial}

\newcommand{\dint}{\dyle\int}

\newcommand{\re}{\mathbb{R}}
\newcommand{\ren}{\re^N}

\newcommand{\dyle}{\displaystyle}

\newcommand{\io}{\int_\O}

\renewcommand{\a }{\alpha }
\renewcommand{\b }{\beta }
\renewcommand{\d }{\rho }
\newcommand{\D }{\Delta }
\newcommand{\e }{\varepsilon }

\newcommand{\g }{\gamma}

\renewcommand{\l }{\lambda }
\renewcommand{\L }{\Lambda }

\newcommand{\n }{\nabla }

\newcommand{\s }{\sigma }

\renewcommand{\O }{\Omega }

\title{Fractional elliptic reaction-diffusion systems with coupled
    gradient terms and different diffusion 
}
\author[1]{Somia Atmani}
\author[2]{Kheireddine Biroud}
\author[3]{Maha Daoud}
\author[4]{El-Haj Laamri \thanks{Corresponding author: el-haj.laamri@univ-lorraine.fr}}

\affil[1]{Laboratoire d’Analyse Non Linéaire et Mathématiques Appliquées, Département de
    Mathématiques, Université Abou Bakr Belkaïd, 13000 Tlemcen, Algeria}

\affil[2]{Laboratoire d’Analyse Non Linéaire et Mathématiques Appliquées, École Supérieure de
    Management de Tlemcen, 13000 Tlemcen, Algeria}
    \affil[3]{
    	POEMS, ENSTA Paris, Institut Polytechnique
    	de Paris, 91120 Palaiseau, France}
\affil[4]{Institut Elie Cartan de Lorraine, Universit\'e de
	Lorraine, 54 506 Vandoeuvre-les-Nancy, France

    \medbreak

    {\normalsize somiaatmani@gmail.com, kh\_biroud@yahoo.fr, maha.daoud@ensta-paris.fr, el-haj.laamri@univ-lorraine.fr}}
\date{ \today}

\begin{document}
    \maketitle

    \begin{abstract}
        In this work, we study the existence and nonexistence of  nonnegative solutions to a class of nonlocal elliptic {systems} set in a bounded open subset of $\R^N$. The diffusion operators are  of type $u_i\mapsto d_i(-\Delta)^{s_i}u_i$ where $0<s_1\neq s_2<1$,  and  the gradients
        of the unknowns act as source terms (see  $(S)$  below). Existence results are obtained by proving  some fine estimates when data belong  to weighted Lebesgue spaces. Those estimates  are new and interesting in themselves.
       
        \bigbreak

        \noindent
        \textbf{Mathematics Subject Classification (2020):}  35R11,  35G50, 35B45, 35J62, 47H10.
        \medbreak
        \noindent
        \textbf{Keywords.} Reaction-diffusion system, Fractional diffusion, Different diffusion, Nonlinear gradient terms, Nonlinear elliptic system,
        Fixed-point theorem, A priori estimates.
    \end{abstract}

    \section{Introduction}
    The main objective of this paper is to investigate the existence and nonexistence of weak solutions {(see Definition \ref{eq:def:faible})} to fractional elliptic systems of type

  	$$ (S)\hspace*{.25cm}
  	        \left\{
        \begin{array}{rclll}
         (-\Delta)^{s_1} u &= & |\nabla v|^{q}+\l f & \text{ in }&\Omega , \\
            (-\Delta)^{s_2} v &= & |\nabla u|^{p}+\mu g & \text{ in }&\Omega , \\
            u =v&=&0 &\hbox{  in }& \mathbb{R}^N\setminus\Omega,
        \end{array}%
        \right.
        $$
    where 
    $0<s_1\neq s_2<1$, $\Omega$ is a bounded regular open subset of $\mathbb{R}^N$ with$N\geq 2$, $p,q\ge 1$, $\lambda,\mu>0$, $f$ and $g$ are  measurable nonnegative functions.
    \\
    Here, the operator $(-\Delta)^s$ denotes the fractional Laplacian introduced par M. Riesz in \cite{Riesz}  and  defined for any $u\in \mathcal{S}(\R^N)$ by:
    
    \begin{equation}\label{eq:fraccionarioel2}
        (-\Delta)^{s}u(x):=a_{N,s}\mbox{ P.V. }\int_{\mathbb{R}^{N}}{\frac{u(x)-u(y)}{|x-y|^{N+2s}}\, dy},\; s\in(0,1),
    \end{equation}
    where P.V. stands for the Cauchy principle value and
    $$a_{N,s}:=\frac{s 2^{2 s} \Gamma\left(\frac{N+2 s}{2}\right)}{\pi^{\frac{N}{2}} \Gamma(1-s)}$$
    is a normalization constant chosen
    so that 
    the following pair of identities :
    $$\lim\limits_{s\to 0^+} (-\Delta)^su=u\quad\text{and }\quad \lim\limits_{s\to 1^-} (-\Delta)^su=-\Delta u$$
    holds. 
    For further details, see, {\it e.g.}, \cite[Proposition 4.4]{dine}
    or \cite[Proposition 2.1]{DaouLaam}.
   
   \medbreak
    
    In the realm of fractional Laplacian operators, two main versions  are often discussed:
    	the so-called regional fractional Laplacian, which we consider in this work (see \eqref{eq:fraccionarioel2}), and the so-called spectral fractional Laplacian defined by:
    		\begin{equation}
    		(-\Delta)^s_S u=\sum_{k \in \mathbb{N}}  \langle u, e_k\rangle \lambda_k^s e_k,
    		\label{eq:spectral_laplacian}
    	\end{equation}
    	where $\lambda_k$ and $e_k$ represent the eigenvalues and eigenfunctions of the classical Laplacian with the homogeneous Dirichlet boundary condition, and $\|e_k\|_{L^2(\Omega)}=1$. 
    	Although these two fractional Laplacian operators are sometimes conflated in the literature, they are fundamentally distinct, as proven in \cite{SerVal} (see also \cite[Section 5.3]{MolRadRaf} , \cite[Remark 5.1]{Molica-et-et-alJMPA}).
    	Moreover, the third and fourth authors have exhibited in \cite[Section 4]{DaouLaam} several simple yet significant examples in dimension 1  that  effectively  highlight the difference between the two operators.
    	     	 
  \noindent   In this context, the same authors and their co-worker previously investigated a class of parabolic reaction-diffusion systems with control of total mass governed by the regional fractional Laplacian, see \cite{DaouLaamBaal2024}.  Subsequently, the third author extended these results to systems involving the spectral fractional Laplacian, as detailed in \cite{DAOUD}. 
  All in all, we have observed that the two operators determine substantial differences in the properties of the systems they govern: they require different approaches, and different technical considerations must be carefully addressed.
Thus, in this paper, we will focus on System~$(S)$ with the regional fractional Laplacian. The same system with the spectral fractional Laplacian will be addressed in a forthcoming work.

    \medbreak 
    More broadly, fractional Laplacians are the
    prototypic nonlocal operators. They occur naturally in diverse
    scientific disciplines, including Biology, Ecology, Finance,
    Materials Science, Nonlinear Dynamics and Anomalous Diffusion. In
    fact, such operators play a fundamental role in modelling complex
    systems and phenomena with long-range interactions. For more
    details, we refer the reader, for instance, to \cite{AbaVal,
        DaouLaam, BucVal}  and references therein.
    \\  Besides, it is worth noting that many types of elliptic systems with gradient terms arise when studying electrochemical models in engineering and some fluid dynamics models. More details can be found in \cite{CRST,Diez} and their references.
    \\
    As for models incorporating fractional diffusion, we refer to
    \cite{MW, KE, KT} for a rigorous justification for introducing the fractional Laplacian into the continuum equation
    of the growing surface as another relaxation mechanism.

    \smallbreak In the present work, we look for \textit{natural}
    conditions on  $(f,\,g)$, $(\l,\,\mu)$ and $(p,\, q)$ which allow
    us to establish {the existence }of solutions to System
    $(S)$.
     \medbreak
 
    Before going further, let us recall some previous results related to our study.
    \\
    $\bullet$ \textbf{Case of a single equation.} System
    $(S)$ is then reduced to
    \begin{equation}\label{eq:Eqel2}
        \left\{
        \begin{array}{rclll}
            (-\Delta)^s u &= & |\nabla u|^{p}+\l f & \text{ in }&\Omega , \\ u &=& 0 &\hbox{  in }& \mathbb{R}^N\setminus\Omega.
        \end{array}%
        \right.
    \end{equation}
    \begin{itemize}
        \item[---] {\it Local case} $s=1$. When $f=0$,  Problem \eqref{eq:Eqel2} is  no other than the stationary 
        viscous {\it Hamilton-Jacobi} Problem (also known as the stationary diffusive {\it Hamilton-Jacobi } Problem), see, {\it e.g.}, \cite{Tran} and references cited therein. When $f\neq0$, we find the stationary generalized {\it Kardar-Parisi-Zhang} {\it(KPZ)} Problem  (see, for instance, \cite{KPZ} and its references). Problem \eqref{eq:Eqel2} has been extensively investigated, see \cite{Ala-Pier, BMP0, BMP1, BMP2, BMP3,BMP4, BGO, GMP, HMV, ADPE,FM1, FM2} and their references. Needless to say that the references provided are not exhaustive, but rather a selection of some relevant sources.
        \item[---] {\it Nonlocal case} $\frac12<s<1$. Although Problem \eqref{eq:Eqel2} has received significant attention when $s=1$, there has been relatively limited research in the case $\frac12<s<1$. In \cite{CV2}, the authors established the existence of solutions when $p<\dfrac{N}{N-2s+1}$. Additionally, they treated the case where $f$ is a Radon measure. In the same paper, as well as in \cite{AB01}, the authors have {considered} natural conditions on $f$ in order to examine the existence of solutions to Problem~\eqref{eq:Eqel2} depending on the value of $p$.
    \end{itemize}
    $\bullet$ \textbf{Case of systems} when $s_1=s_2=s$.
    \begin{itemize}
        \item[---] {\it Local case} $s=1$. Unlike the case of single equation, few results are known {for} systems where the gradients of the unknowns act as source terms. Let us review the known results {regarding} such systems.  In~\cite{AttaBentLaam}, the fourth author and his coworkers have established
        the existence of nonnegative solutions to System $(S)$  for any $(p,q)\in [1,+\infty)^2$ such that $pq>1$, under  suitable integrability conditions on $f$ and $g$,  and smallness conditions on $\lambda$ and $\mu$. In addition, nonexistence results were also identified  for $\lambda$ (or  $\mu$) large, or $f$ and $g$ of low integrability. Similar results was obtained  in \cite{AAL11} about an elliptic {system} where the right-hand sides are potential and gradient terms, namely, $-\Delta u = v^{q}+\lambda f \;,\; -\Delta v =  |\nabla u|^{p}+\mu g$, as well as its parabolic version. Moreover, we refer to \cite{AtB} where the authors studied the particular case  $0<p\leq 2$, $q>0$ with $pq<1$, and $\lambda=\mu=1$.
        Many other studies, such as  \cite{BOPorr,BOPu,Diez}, have also investigated the existence of solutions to local elliptic systems with gradient terms.
        \item[---] {\it Nonlocal case} $\frac12<s_1=s_2=s<1$.
        In this case, we have 
         studied System~$(S)$ in \cite{AtmaBirDaouLaam}, and our work appears to be a novel contribution. For the reader's convenience, let us summarize the main existence and nonexistence results that we obtained.
        \\
        $\diamond$ The existence of solutions was obtained  for
        $(p,q)\in \Big[1,\dfrac{s}{1-s}\Big)^2$, under appropriate integrability assumptions on $f$ and $g$ and smallness conditions on $\lambda$ and $\mu$.\\
        $\diamond$ Nonexistence results was  proven in the following three cases :
        \begin{itemize}
            \item[1)] $\lambda$ and $\mu$ are large even if $1<p,q< \dfrac{s}{1-s}$;
            \item[{2)}] $\max\{p,q\}\geq\dfrac{1}{1-s}$ even if $f$ and  $g$ are bounded and  no matter how small are $\lambda$ and $\mu$;
            \item[3)] $f$ and $g$ have less regularity.
        \end{itemize}
  In the same work \cite[Section 5]{AtmaBirDaouLaam}, we have addressed other nonlinear systems
    	with different structures. Notably, we have treated systems with potential-gradient reaction terms of the form		$(-\mathrm{\Delta})^{s} u =  v^{q}+\lambda    f\;,\; (-\mathrm{\Delta})^{s} v =  |\nabla u|^{p}+\mu g.$ 
        Moreover, it is noteworthy that we have recently  investigated in~\cite{AtmaBirDaouLaam-FCAA} the parabolic version  of System~$(S)$.
 \end{itemize}
    
%

     
  \smallbreak  
    
  
 	
 	Systems of type~$(S)$ can be seen as a generalization of the nonlocal $KPZ$ stationary equation. 
	As in the latter model, it seems to be natural to consider $s_1, s_2 > \frac{1}{2}$; see \cite[Section 1]{BJ} for a detailed justification of this condition, especially in the presence of the gradient acting as a source term.
   \smallbreak  	
 From now on, unless explicitly stated otherwise,  we shall assume $\frac12<s_1 < s_2<1$. The reverse case $s_2 < s_1$ could be handled analogously, as $(s_1,p)$ and $(s_2,q)$ play symmetrical roles.
 	To the best of our knowledge, System~$(S)$ under this assumption has not been investigated yet.   
	{This is the main goal of this paper},
 	which complements our earlier study in~\cite{AtmaBirDaouLaam}. On the other hand, we have recently examined in~\cite{AtmaBirDaouLaam-Parab-s1-s2} the parabolic version of System~$(S)$, extending our previous work~\cite{AtmaBirDaouLaam-FCAA}. 
 In both the elliptic and the parabolic cases,  the extension from $s_1=s_2$ to $s_1\ne s_2$ turns out to be much more difficult than it appears at first sight. Indeed, in addition to the difficulties faced in the
    case $s_1=s_2$, there are other challenges that arise due to the fact that the diffusion operators are different.
Hence, the techniques used in the local case are no longer
applicable to handle the new situation. In fact, the major
difficulties are coming from the loss of regularity caused by the
{fractional Laplacian} and the lack of the homogeneity of this
operator. In order to overcome these {difficulties,} we will use
in a convenient way (as we have done in \cite{AtmaBirDaouLaam}) :
\\
(i) sharp estimates on Green's function of the fractional Laplacian obtained in~\cite{BJ};\\
(ii) fine weighted estimates on the solution to the following fractional Poisson problem
$$ (-\Delta)^{s_2} w =  h  \mbox{ \mbox{ in}  } \O, \quad w =  0  \mbox{  in }  \mathbb{R}^N\backslash\O,$$
 where $h$  is  a given  function such that $h\d^{a(1-s_1)} \in L^m(\O)$ where $m\ge 1$, $0\le a<\dfrac{s_2}{1-s_1}$ and $\rho(x): =\text{dist}(x,\partial \Omega)$ ; see Theorems \ref{eq:regularity_usurdeltas} and  \ref{eq:main_regularity_resultel2} below.\\
It should be  mentioned that these estimates  are new, have their own
interest and can serve to address other similar nonlocal problems. They also improve those obtained in~\cite{AtmaBirDaouLaam, CV1,CV2,AB01} and extend those of \cite{BV} in the local case. 
{ This is a further essential goal of this work.}
 
\medbreak
\noindent Relying on these estimates, the existence of {a solution} will be  obtained for $1\leq q<\dfrac{s_1}{1-s_1}$ and
$1\leq p<\dfrac{s_2}{1-s_1}$ under suitable regularity assumptions
on $f$ and $g$, and smallness conditions on $\l $ and $\mu $, see
Theorems \ref{eq:existence_theorem_for_Sel2} and \ref{eq:pq=1}.
Nonexistence results will be also obtained in these three cases :
\begin{itemize}
    \item[{(i)}] $\lambda$ and $\mu$ are large even if $1< q<\dfrac{s_1}{1-s_1}$ and $1< p<\dfrac{s_2}{1-s_1}$, see Theorem~\ref{eq:NonexisThm1};
    \item[{(ii)}]  $\max\{p,q\}\geq \dfrac{1}{1-s_1}$ even if $f$ and  $g$ are bounded and  no matter how
    small are $\lambda$ and $\mu$, see Theorem~\ref{eq:NonExisThm2-};
    \item [{(iii)}] $f$ and $g$ are less regular, see Theorem \ref{eq:Nonexistence3}.
\end{itemize}
{Let us} point out that our nonexistence  Theorems~\ref{eq:NonexisThm1}, \ref{eq:NonExisThm2-} and \ref{eq:Nonexistence3} show, in some sense,
that our assumptions on the data in Theorem~\ref{eq:existence_theorem_for_Sel2} are indeed necessary.
\smallbreak
\noindent Regarding the existence, we believe that  Theorem~\ref{eq:existence_theorem_for_Sel2}  holds for $q\in \Big[\dfrac{s_1}{1-s_1},\dfrac{1}{1-s_1}\Big)$ and $p\in \Big[\dfrac{s_2}{1-s_1},\dfrac{1}{1-s_1}\Big)$. To date, we have not been able to prove the existence in this case. We leave it as an open problem.
As a matter of fact, this question is still open even in the case
$s_1=s_2$ (see~\cite{AtmaBirDaouLaam}),  as well as in the case of a
single equation (see \cite{AB01}).

\medbreak


\medbreak

The rest of this paper is organized as follows. In Section 2, we
give some useful tools such as {regularity} results for nonlocal
Poisson problem with general datum in $L^1$ or in the space of
Radon measures. Furthermore, we recall important estimates on
Green's function, that will be used to prove our existence theorems.
In Section 3, we prove some new regularity and compactness results
in weighted spaces. In Section  4, we prove the main existence
results of nonnegative solutions to System~$(S)$
in suitable Sobolev spaces. Moreover, we present some nonexistence
results that show, in some sense, the optimality of the
assumptions considered in our existence theorems. 
\\

\noindent{\bf Notations.}\quad Let us fix some notations which will be used throughout this paper.
\begin{itemize}
    %
    %
      \item[---] For any real numbers $\varsigma \geq 0$ and  $k>0$,  let  $T_k(\varsigma):
	=\max\{-k, \min\{k,\varsigma\}\}$.  
    \item[---] For $\beta>0$,  $L^1(\Omega, \rho^\beta(x)dx):=\{h:\Omega\rightarrow \mathbb{R} \text{ measurable }\;;\; {\|h\|_{L^1(\Omega, \rho^\beta)}}<+\infty\}$, {where} $$\|h\|_{L^1(\Omega, \rho^\beta)}=\int_{\Omega} |h(x)|\rho^\beta(x)\,dx.$$
    \item[---] Let $\varphi$ and $\psi$ be two nonnegative functions.  We
    write $\varphi\simeq \psi$ if there exist two positive constants
    $C_1$ and $C_2$ such that
    $$C_1\varphi\leq \psi\leq C_2\varphi.$$
    	\item[---] {By $C$, we denote a positive constant which may be different from line to another.}
    	\item[---] For any $p>1$, we denote by $p'$ its conjugate, that is $\frac{1}{p}+\frac{1}{p'}=1.$
    	\item[---] Let $x_0\in\mathbb{R}^N$ and $r>0$. $B_r(x_0)$ denotes the open ball centered at $x_0$ with radius $r$.
\end{itemize}
\hfill$\square$

\section{Preliminaries and tools}
To make this paper self-contained, we review in this section some classical results that will be used later on.

\subsection{
{Sobolev Spaces}}

{For the reader’s convenience, we recall here some classical definitions}. Let $s\in (0,1)$ and $p\in [1,+\infty)$.

\noindent$\bullet$ { \underline{Fractional Sobolev {space} $W^{s,p}(\mathbb{R}^N)$}}\\
The space
$$
W^{s,p}(\mathbb{R}^N):= \left\{ \phi \in L^p(\mathbb{R}^N) \; ; \;
\iint_{\mathbb{R}^N \times \mathbb{R}^N}
\frac{|\phi(x)-\phi(y)|^p}{|x-y|^{N+ps}}dx dy  < +\infty \right\}, $$
endowed with  its canonical norm
$$
\|\phi\|_{W^{s,p}(\mathbb{R}^N)} := \left( \|\phi\|_{L^p(\mathbb{R}^N)}^p + \iint_{\mathbb{R}^N \times \mathbb{R}^N} \frac{|\phi(x)-\phi(y)|^p}{|x-y|^{N+ps}} dx dy  \right)^{\frac{1}{p}}\;
$$
is a Banach space.

\medbreak

\noindent$\bullet$ { \underline{Fractional  Sobolev {space} ${\Bbb W}_0^{s,p}(\Omega)$}}\\
The space $\mathbb{W}_0^{s,p}(\O)$ is defined by
$$
\mathbb{W}_0^{s,p}(\O):= \left\{\phi\in W^{s,p}(\mathbb{R}^N)\;; \phi= 0 \textup{ in } \mathbb{R}^N\setminus \Omega\right\}.
$$
Endowed with the norm induced by $\|\cdot\|_{W^{s,p}(\R^N)}$ which we denote $|\cdot|_{\mathbb{W}^{s,p}_0(\Omega)}$,  $\mathbb{W}_0^{s,p}(\O)$ is a Banach space.

\medbreak

\noindent For $\phi\in \mathbb{W}_0^{s,p}(\O)$, we set
\[ \|\phi\|_{\mathbb{W}^{s,p}_0(\O)} := \left( \iint_{D_{\O}} \frac{|\phi(x)-\phi(y)|^p}{|x-y|^{N+sp}} dx dy  \right)^{1/p},\]
where
\[ D_{\O} := (\mathbb{R}^N \times \mathbb{R}^N) \setminus \left((\R^N\setminus \O) \times (\R^N\setminus \O)\right)  .\]

\noindent It is clear that $\|\cdot\|_{\mathbb{W}^{s,p}_0(\O)}$  defines a norm on $\mathbb{W}_0^{s,p}(\O)$. In addition, $\O$ is assumed  to be bounded, then  the norms $|\cdot|_{\mathbb{W}_0^{s,p}(\O)}$ and $\|\cdot\|_{\mathbb{W}_0^{s,p}(\O)}$ are equivalent.
{\medbreak
\noindent$\bullet$ {\underline{Sobolev spaces ${W}^{1,p}(\mathbb{R}^N)$ and ${\Bbb W}_0^{1,p}(\Omega)$}}\\
The space $W^{1,p}(\mathbb{R}^N)$ is defined as follows :
$$
W^{1, p}(\mathbb{R}^N):=\left\{\phi :\mathbb{R}^N\to\R \;\text{measurable}\;;\|\phi\|_{W^{1, p}(\mathbb{R}^N)}<+\infty   \right\},
$$
where
$$
\|\phi\|_{W^{1, p}(\mathbb{R}^N)}:=\left(\sum_{|\alpha| \leq 1} \int_{\mathbb{R}^N}\left|D^\alpha \phi(x)\right|^p  dx \right)^{1 / p}.
$$
By definition, we have
\[ \mathbb{W}_0^{1,p}(\Omega):= \left\{\phi\in W^{1,p}(\mathbb{R}^N)\;;\; \phi= 0 \textup{ in } \mathbb{R}^N\setminus \O\right\}.\square \]
}

\subsection{ Some properties of Green's function}
For a given $x\in \O$,  the Green's function  $\mathcal{G}_s(x ,y)$ is the solution to the following problem:
\begin{equation}\label{eq:Green0el2}
    \left\{\begin{array}{rclll} (-\Delta)^s_y \mathcal{G}_s(x ,y) &= & \delta_x & \mbox{  if  }&y\in \O,\\ \mathcal{G}_s(x ,y) &= & 0 & \mbox{ if }& y\in
        \mathbb{R}^N\backslash\O,
    \end{array}
    \right.
\end{equation}
where {$\delta_x$} is the Dirac mass at $x$.

\medbreak

Here are some useful properties of this
function and its gradient.  We  refer  the interested readers to
\cite{AB01,CZ,BJ0,BJ} for a detailed  proof (which relies on
probabilistic tools).
\begin{lem}\label{eq:Green_propertiesel2} 
    Assume that  $s\in (0,1)$. Then, for \textit{a.e.} $x ,y\in \O$, we have
    \begin{equation}\label{eq:First_GPel2}
        \mathcal{G}_s(x ,y)\simeq
        \frac{1}{|x-y|^{N-2s}}\bigg(\frac{\d^s(x)}{|x-y|^{s}}\wedge
        1\bigg) \bigg(\frac{\d^s(y)}{|x-y|^{s}}\wedge 1\bigg).
    \end{equation}
    In particular, there exists $C>0$ such that
    \begin{equation}\label{eq:Second_GPel2}
        \mathcal{G}_s(x ,y)\le C\min\left\{\frac{1}{|x-y|^{N-2s}},
        \frac{\d^s(x)}{|x-y|^{N-s}}, \frac{\d^s(y)}{|x-y|^{N-s}}\right\}.
    \end{equation}
    Moreover, for any $\eta\in (0,1)$, we get
    \begin{equation}\label{eq:Third_GPel2}
        \mathcal{G}_s(x ,y)\le C\frac{\d^{\eta s}(y)\d^{(1-\eta)
                s}(x)}{|x-y|^{N-s}}\mbox{  and   } \hspace{0.2cm}\mathcal{G}_s(x ,y)\le
        C\frac{\d^{\eta s}(y)}{|x-y|^{N-s(2-\eta)}}.
    \end{equation}
    If $s\in (\frac 12,1)$, then there exists $C>0$ such that
    \begin{equation}\label{eq:Fourth_GPel2}
        |\n_x \mathcal{G}_s(x ,y)|\le C
        \mathcal{G}_s(x ,y)\max\left\{\frac{1}{|x-y|}, \frac{1}{\d(x)} \right\}.
    \end{equation}
    Therefore,
    \begin{equation}\label{eq:Fifth_GPel2}
        |\n_x \mathcal{G}_s(x ,y)|\le
        \frac{C}{|x-y|^{N-2s+1}\d^{1-s}(x)}.\square
    \end{equation}
\end{lem}

\subsection{ Fractional Poisson equation:  regularity and compactness}
Let us consider  the {following problem}
{\begin{equation}\label{eq:FPEel2}
        \left\{\begin{array}{rclll}
            (-\Delta )^s w&=&h &\hbox{   in   }&  \Omega,\\ w&=&0   &\hbox{   in   }& \mathbb{R}^N\setminus\Omega,
        \end{array}\right.
\end{equation}}
where $0<s<1$ and $h$ is a measurable function such that $h\in L^m(\O)$ with $m\geq 1$.
\smallbreak
\noindent $\bullet$ {\bf Weak solution.} Let us proceed by specifying what we mean by a (weak) solution to  Problem~\eqref{eq:FPEel2}.
\begin{defi}\label{eq:def1o}  
    We say that $w$ is a (weak) solution to Problem \eqref{eq:FPEel2} if
    $w\in L^1(\O)$, and for any $\phi\in \mathbb{X}_s(\Omega)$, we have
    $$
    \io w(-\Delta )^s\phi\,  dx =\io \phi h\, dx,
    $$
    where
    $$
    \mathbb{X}_s(\Omega) := \Big\{\phi\in
    \mathcal{C}(\mathbb{R}^N)\,;\,\text{supp}(\phi)\subset \overline{\O},\, (-\Delta )^s\phi\, \text{ exists} \text{
        and } (-\Delta )^s\phi\in L^\infty(\O)\Big\}.
    $$
\end{defi}

\noindent $\bullet$ {\bf Green's function representation.}
The solution defined above is given by the following  representation formula:
\begin{equation}\label{eq:representation_formulael2}
	w(x)=\int_\Omega    \mathcal{G}_s(x ,y)h(y) dy.
\end{equation}

The following theorem summarizes the existence result and some essential properties of these solutions.
\begin{thm}[{\cite[Theorem 23]{LPPS}}]\label{eq:Existence_entropy_solutionel2} 
Assume that $h\in L^1(\O)$.
\begin{enumerate}
 \item[\em{1.}] Problem~\eqref{eq:FPEel2}  has a unique weak  solution $w$.
 \item[\em{2.}] This solution  can be constructed in the following way :\\
	For any $n\in \mathbb{N}^*$, let $w_n$ be the 
	solution to the approximating problem
	\begin{equation}\label{eq:Approximating_Proel2}
		\left\{\begin{array}{rclll} (-\Delta)^s w_n &= & T_n(h) & \mbox{\em
				in  } &\O,\\ w_n &= & 0 & \mbox{\em in } & \mathbb{R}^N\backslash\O;
		\end{array}
		\right.
	\end{equation}
it exists as the integral obtained by replacing $h$ with $T_n(h)$ in~\eqref{eq:representation_formulael2}.\\
The sequence $(w_n)_n$ converges in the following sense:
	\begin{equation}\label{eq:L1duuel2}
		w_n\to w\mbox{  strongly in } \mathbb{W}^{s,\g}_0(\O), \qquad \forall\g<\frac{N}{N-2s+1}
	\end{equation}
and
	\begin{equation} \label{eq:tkuel2}
	\forall k>0,\quad	T_k(w_n)\to T_k(w)\hbox{  strongly in }
		\mathbb{W}^{s,2}_{0}(\Omega).
	\end{equation}
 \item[\em{3.}]Moreover,
	\begin{equation}  \label{eq:L1uel2}
		w\in L^\theta(\Omega) \,, \qquad  \forall  \ \theta\in \big(1,
		\frac{N}{N-2s}\big)
	\end{equation}
	and
	\begin{equation}\label{eq:L1duel2}
		(-\Delta)^{\frac{s}{2}}   w \in L^r(\Omega) \,, \qquad
		\forall  \  r \in \big(1,  \frac{N}{N-s} \big) \,.
	\end{equation}

\end{enumerate}
		\end{thm}

\medbreak

There also holds several regularity results in weighted spaces.
\begin{prop}[{\cite[Proposition 2.3]{CV2}}]\label{eq:Existence_with_weightel2} 
    Let $s\in (\frac 12,1)$. Assume that $h\in L^1(\Omega,\d^\beta)$ with $0\leq \beta<2s-1$.
    Then, Problem \eqref{eq:FPEel2} has a unique weak solution $w \in
    {\mathbb{W}_0^{2s-\gamma,p}(\O)}$ for any {$p\in(1,
        \frac{N}{N+\beta-2s})$} where $\gamma=\beta+\frac
    N{p'}$ if $\beta>0$ and $\gamma>\frac N{p'}$ if $\beta=0$.\\
    Moreover, there exists $C:=C(p,s,\O)>0$ independent of $w$ and
    $h$ such that
    $$\|w\|_{\mathbb{W}_0^{2s-\gamma,p}(\O)}\le C\|h\|_{L^1(\Omega,\d^\beta)}.$$
\end{prop}
Also, for a datum $h\d^{a(1-s)}\in L^m(\O)$ with $a\ge 0$ and $m\ge 1$,  we have:
\begin{thm}[{\cite[Theorem 3.2]{AtmaBirDaouLaam}}] \label{eq:regubeta1}
    Let  $s\in (\frac 12,1)$ and  $h$ be a measurable function such that $h\rho ^{a(1-s)}\in L^{m}(\O)$ with $m\ge 1$ and $0\le a<\dfrac{s}{1-s}$. Let $w$ be  the unique weak solution to Problem~\eqref{eq:FPEel2}.
    Moreover, let us set
    $\widehat{m}_{s,a}:=\frac{mN}{(N-m(s-a(1-s)))_+}$ and $\widetilde{m}_{s,a}=:\min\{\widehat{m}_{s,a},
    \frac{mN}{(N-m(2s-1))_+}\}$. {Then, we have the following assertions. }

    \noindent
    \begin{enumerate}

        \item[\em (i)] If
        $m<\max\{\frac{N}{2s-1},\frac{N}{s-a(1-s)}\}$, then $|\n w|\d^{1-s} \in L^{\g}(\O)$ for any
        $\gamma<\widetilde{m}_{s,a}$ and there exists $C>0$ such that
        \begin{equation}  \label{eq:comp_result101} 
            \left\||\n w|\d^{1-s}\right\|_{L^{\g}(\O)}\leq C\bigg(
            \left\|\dfrac{w}{\d^s}\right\|_{L^{\widehat{m}_{s,a}}(\O )}+
            \|\Im(h\d^{a(1-s)})\|_{L^{\gamma}(\O )}\bigg),
        \end{equation}
        where $\Im: L^m(\O)\to L^\g(\O)$  is the linear operator defined by
        \begin{equation}\label{eq:Pi1}
            \Im(g)(x):=
            \left\{
            \begin{array}{lll}
                \dyle \io\frac{g(y)}{|x-y|^{N-(s-a(1-s))}} dy
                &
                \mbox{\em if } & a\ge 1,\\ \\
                \dyle \io\frac{g(y)}{|x-y|^{N-2s+1}} dy &\mbox{\em
                    if
                } & 0\le a<1.
            \end{array}
            \right.
        \end{equation}
        Moreover for any $\g<\widetilde{m}_{s,a}$, there exists a constant $C>0$ such that

        \begin{equation}\label{eq:mainalpha} 
            \left\||\n w|\d^{1-s}\right\|_{L^{\g}(\O)}\leq C\|h \d^{a(1-s)}\|_{L^{m}(\O )}.
        \end{equation}
        \item[\em (ii)] If
        $m=\max\{\frac{N}{2s-1},\frac{N}{s-a(1-s)}\}$, then $|\n
        w|\d^{1-s}\in L^\g(\O)$ for any $\g<+\infty$ and there exists $C>0$ such that

        \begin{equation}\label{eq:comp_result01221}  
            \left\||\n w|\d^{1-s}\right\|_{L^{\g}(\O)}\leq C\bigg(
            \left\|\dfrac{w}{\d^s}\right\|_{L^{\g}(\O )}\,\,\,
            + \|\Im(h\d^{a(1-s)})\|_{L^{\g}(\O )}\bigg).
        \end{equation}
        \item[\em (iii)] If
        $m>\max\{\frac{N}{2s-1},\frac{N}{s-a(1-s)}\}$, then $|\n
        w|\d^{1-s}\in L^\infty(\O)$ and there exists $C>0$ such that

        \begin{equation}\label{eq:comp_result011}  
            \left\||\n w|\d^{1-s}\right\|_{L^{\infty}(\O)}\leq C\bigg(
            \left\|\dfrac{w}{\d^s}\right\|_{L^{\infty}(\O )}\,\,\,
            + \|\Im(h\d^{a(1-s)})\|_{L^{\infty}(\O )}\bigg).
        \end{equation}

    \end{enumerate}
    Moreover,  for $\g<\min\{\frac{N}{N-(s-a(1-s))},
    \frac{N}{N-(2s-1)}\}$ fixed, the operator
    $\Gamma:L^1(\O,\d^{a(1-s)})\to
    \mathbb{W}^{1,\g}_0(\Omega,\d^{\g(1-s)})$ defined by $\Gamma(h)=w$,
    where $w$ is the unique weak solution to \eqref{eq:FPEel2}, is
    compact.$\square$
\end{thm}

\subsection{ Other miscellaneous useful tools }
\subsubsection{{\bf Stein's Lemma}}
\begin{lem}[{\cite[Chapter V, Section 1.2]{Stein}}]\label{eq:Stein}
    Let $0<\lambda<N$,  $1\le p<\ell<+\infty$ be such that  $\dfrac{1}{\ell}+1=\dfrac{1}{p}+\dfrac{\l}{N}$. For $g\in L^p(\ren)$, we define $$J_\lambda(g)(x):=\int_{\ren}
    \dfrac{g(y)}{|x-y|^\l}dy.$$
    \begin{itemize}
        \item[\em a)] $J_\lambda$ is well defined in the sense that the integral converges absolutely for  almost any $x\in \mathbb{R}^N$. \item[\em b)] If $p>1$, then
        $\| J_\lambda(g)\|_{L^\ell(\mathbb{R}^N)}\le C_{p,\ell}\|g\|_{L^p(\mathbb{R}^N)}$. \item[\em c)] If $p=1$, then $\left|\{x\in \mathbb{R}^N\,;\,| J_\lambda(g)(x)|>\sigma\}\right|\le \bigg(\dfrac{
                    C_{1,\ell}\|g\|_{L^1(\mathbb{R}^N)}}{\sigma}\bigg)^\ell$.
    \end{itemize}
\end{lem}

\subsubsection{\bf Interpolation Inequality}

\begin{prop}[{\cite[Theorem 1]{BL}}]\label{eq:interpolation} 
    Let $h\in L^1(\Omega)$ and $\{h_n\}_n\subset L^r(\Omega)$ with  $1< r<+\infty$. Assume    that $h_n\rightarrow h$ strongly in $L^1(\Omega)$ and $\{h_n\}_n$ is bounded in $L^r(\Omega)$.\\ Then,
    $h\in L^r(\Omega)$, and  for any $a\in [1,r)$  we have the following interpolation inequality
    \begin{equation} \label{eq:interpolationInequalityel2}
        \|h_n-h\|_{L^a(\Omega)}\leq \|h_n-h\|_{L^1(\Omega)}^\theta \|h_n-h\|_{L^r(\Omega)}^{1-\theta} \text{ where }  \theta:= \displaystyle\frac{r-a}{a(r-1)}.
    \end{equation}
    In particular, $h_n\rightarrow h$  strongly in $L^a(\Omega)$.
\end{prop}
\subsubsection{ Weighted Hardy inequality}

\begin{thm}[{\cite[Theorem 1.6]{Necas}} or {\cite[Theorem 3.1]{EdmundsHurri}}]\label{eq:HardyInequalityel2}
    Let $\nu \in (1,+\infty)$ and $0\leq\gamma<\nu-1$. Then, there
    exists a positive constant $C=C(\O, \gamma,\nu)$ such that for any
    $\varphi\in C_0^{\iy}(\O)$, we have
    \begin{equation}\label{eq:HardyIneqel2}
        \io \d^{\gamma-\nu}(x)|\varphi(x)|^\nu  dx \leq C \io \d^{\gamma}(x)|\n \varphi(x)|^{\nu} dx.
    \end{equation}
\end{thm}

\section{Regularity and compactness results in weighted spaces}\label{eq:sec:regu-comp}
In this section, we establish two  weighted estimates on the solution to the following fractional Poisson problem
\begin{equation}\label{eq:poissonel2}
	\left\{\begin{array}{rclll} (-\Delta)^{s_2} w &= & h & \mbox{
			\mbox{in}  } &\O,\\ w &= & 0 & \mbox{ in } & \mathbb{R}^N\backslash\O,
	\end{array}
	\right.
\end{equation}	
where $h\d^{a(1-s_1)}\in L^m(\O)$, $m\geq 1$ and $0<s_1<s_2<1$.
\smallbreak
As mentioned in the introduction, these results are new, have their own interest and can serve to address other similar nonlocal problems.
 The first theorem applies to any $0<s_1<s_2<1$, while the second one applies to $\frac12<s_1<s_2<1$, as it relates to the gradient term. It is noteworthy to recall that the equality case $s_1=s_2$ has already been thoroughly examined in our previous work~\cite{AtmaBirDaouLaam}.  However,  the case $s_1<s_2$ is more demanding and needs some significant changes.
 
\smallbreak

\begin{thm}\label{eq:regularity_usurdeltas} 
 {Let $s_1, s_2\in (0,1)$ such that
   	$s_1<s_2$.} Let $h$  be a measurable function such that $h\d^{a(1-s_1)}\in
    L^m(\O)$, where $m\ge 1$ and $0\le a<\dfrac{s_2}{1-s_1}$.
    Let $w$ be the unique solution to Problem \eqref{eq:poissonel2}. Then,  there exists a positive constant $C$ such that
    \begin{enumerate}
        \item[\em (i)] if $1\leq m<\frac{N}{2s_2-s_1-a(1-s_1)}$,   $\dfrac{w}{\d^{s_1}}\in
        L^\g(\O)$ for any $\g\leq \widehat{m} _{s_1,s_2,a}:=
        \frac{mN}{N-m(2s_2-s_1-a(1-s_1))}$. Moreover,
        \begin{equation}\label{eq:hardyuuel2}
            \left\|\frac{w}{\d^{s_1}}\right\|_{L^{\g}(\O)} \leq
            C\|h\d^{a(1-s_1)}\|_{L^{m}(\O)}\, ;
        \end{equation}
        \item[\em (ii)] if $m=\frac{N}{2s_2-s_1-a(1-s_1)}$,  $\dfrac{w}{\d^{s_1}}\in
        L^\g(\O)$ for any $\g<+\infty$ and
        \begin{equation}\label{eq:hardyuu0el2}
            \left\|\frac{w}{\d^{s_1}}\right\|_{L^{\g}(\O)} \leq
            C\|h\d^{a(1-s_1)}\|_{L^{m}(\O)}\; ;
        \end{equation}
        \item[\em (iii)] if $m>\frac{N}{2s_2-s_1-a(1-s_1)}$, $\dfrac{w}{\d^{s_1}}\in
        L^\infty(\O)$ and
        \begin{equation}\label{eq:hardyuu022el2}
            \left\|\frac{w}{\d^{s_1}}\right\|_{L^{\infty}(\O)} \leq
            C\|h\d^{a(1-s_1)}\|_{L^{m}(\O)}.
        \end{equation}
    \end{enumerate}
\end{thm}

\noindent\begin{proof}
    We give the proof for the  case~(i). The demonstration  of the other cases is quite similar. We closely follow the arguments used in \cite{AtmaBirDaouLaam}. 

{Without loss of generality,}    assume that $h\geq 0$.
    Moreover, let  $w$ be the unique solution to Problem~\eqref{eq:poissonel2}. Then, from the representation formula, we get
    $$ w(x)=\io \mathcal{G}_{s_2}(x ,y) h(y)\,dy. $$
    Thus,
    \begin{equation}\label{eq:hardyyel2}
        \left\{
        \begin{array}{lll}
            \left|\dfrac{w(x)}{\d^{s_1}(x)}\right| & \le  & \dyle \io
            \frac{\mathcal{G}_{s_2}(x ,y)}{\d^{s_1}(x)}h(y)\,dy,\\
            & \le& \dyle\int_{\{z\in\O\;;\;|x-z|<\frac 12 \d(x)\}}
            \frac{\mathcal{G}_{s_2}(x ,y)}{\d^{s_1}(x)}h(y)\,dy + \int_{\{z\in\O\;;\;\frac 12
                \d(x)\le |x-z|<\d(x)\}}
            \frac{\mathcal{G}_{s_2}(x ,y)}{\d^{s_1}(x)}h(y)\,dy
            \\ \\&+ &  \dyle\int_{\{z\in\O\;;\;\d(x)\le |x-z|<
                \d(z)\}} \dfrac{\mathcal{G}_{s_2}(x ,y)}{\d^{s_1}(x)} h(y) dy\\ \\
            & +&\dyle \int_{\{z\in\O\;;\;|x-z|\ge \max\{\d(z),\d(x)\}\}} \dfrac{\mathcal{G}_{s_2}(x ,y)}{\d^{s_1}(x)}h(y) dy,\\\\
            &=& J_1(x)+J_2(x)+J_3(x)+J_4(x).
        \end{array}
        \right.
    \end{equation}
    $\bullet$ Let us start by estimating $J_1$. Notice that, $\d$ is a Lipschitz function. Then,  for any $x \in\O$ and $y\in \{z\in\O\;;\;|x-z|<\frac 12
    \d(x)\}$, we have $|\d(x)-\d(y)|\le |x-y|\le \frac 12 \d(x).$
    Hence,
    $$
    |x-y|\le \frac 12 \d(x)\le \d(y)\le \frac 32 \d(x).
    $$
    Furthermore,  from  \eqref{eq:First_GPel2} we get  $\mathcal{G}_{s_2}(x ,y)\le \dfrac{C}{|x-y|^{N-2s_2}}$.
    Thus,
    \begin{equation}
        \left\{\begin{array}{llll}
            J_1(x) & \leq & C\dyle\int_{\{z\in\O\;;\;|x-z|<\frac 12 \d(x)\}}
            \frac{h(y)}{\d^{s_1}(x)|x-y|^{N-2s_2}}\,dy, \vspace{0.2cm}\\
            & \leq & C\dyle\int_{\{z\in\O\;;\;|x-z|<\frac 12 \d(x)\}}
            \frac{h(y)}{|x-y|^{N-2s_2+s_1}}\,dy,\vspace{0.2cm}\\
            &\le & C\dyle \int_{\{z\in\O\;;\;|x-z|<\frac 12 \d(x)\}}
            \frac{h(y)\d^{a(1-s_1)}(y)
                \d^{-a(1-s_1)}(y)}{|x-y|^{N-2s_2+s_1}}\,dy,\vspace{0.2cm}\\
            &\le & C\dyle \int_{\{z\in\O\;;\;|x-z|<\frac 12 \d(x)\}}
            \frac{h(y)\d^{a(1-s_1)}(y)
            }{|x-y|^{N-(2s_2-s_1-a(1-s_1))}}\,dy,\vspace{0.2cm}\\
            &\le & C\dyle \int_{\O}
            \frac{h(y)\d^{a(1-s_1)}(y)
            }{|x-y|^{N-(2s_2-s_1-a(1-s_1))}}\,dy.
        \end{array}
        \right.
    \end{equation}
    Since $h\d^{a(1-s_1)}\in L^m (\O)$, then by using Lemma \ref{eq:Stein}, it follows that $J_{1}\in L^{ \widehat{m} _{s_1,s_2,a}}(\O)$ and
    $$
    \|J_{1}\|_{L^{\g}(\O)}\le C\|h\d^{a(1-s_1)}\|_{L^{m}(\O)}\qquad\hbox{
        for any } \g\leq  \widehat{m} _{s_1,s_2,a}.
    $$
    $\bullet$ Now, we  deal with $J_2$. Using estimate \eqref{eq:Third_GPel2}
    with $\eta=a\frac{1-s_1}{s_2}<1$, it holds that
    $$
    \mathcal{G}_{s_2}(x ,y)\le
    C\frac{\d^{a(1-s_1)}(y)\d^{s_2-(1-s_1)a}(x)}{|x-y|^{N-s_2}}.
    $$
    Therefore, it follows that
    \begin{equation}
        \left\{\begin{array}{llll}
            J_2(x) &= & \dyle \int_{\{z\in\O\;;\;\frac 12 \d(x)\leq|x-z|<\d(x)\}}\frac{\mathcal{G}_{s_2}(x ,y)}{\d^{s_1}(x)}h(y)\,dy,\vspace{0.2cm}\\
            &\le &\dyle  C\int_{\{z\in\O\;;\;\frac 12 \d(x)\leq |x-z|<\d(x)\}}
            \frac{h(y)\d^{a(1-s_1)}(y)\d^{s_2-a(1-s_1)}(x)}{|x-y|^{N-s_2+s_1}}\,dy,\vspace{0.2cm}
            \\
            &\leq& C \dyle\int_{\{z\in\O\;;\;\frac 12 \d(x)\leq|x-z|<\d(x)\}} \frac{h(y)\d^{a(1-s_1)}(y)}{|x-y|^{N-(2s_2-s_1-a(1-s_1))}}\,dy,\vspace{0.2cm}\\
            &\leq& C \dyle\int_{\O} \frac{h(y)\d^{a(1-s_1)}(y)}{|x-y|^{N-(2s_2-s_1-a(1-s_1))}}\,dy.
        \end{array}
        \right.
    \end{equation}
    Then, we conclude {as for} the term $J_1$.\\
    \noindent $\bullet$ {Let us now} treat $J_3$. Since $s_2>s_1$, then by
    using estimate \eqref{eq:First_GPel2} we get
    \begin{equation}\label{eq:J31(x)el2}
        \left\{\begin{array}{ll}
            J_{3}(x) &\displaystyle \le  C\frac{1}{\d^{s_1}(x)} \int_{\{z\in\O\;;\;\d(x)\le |x-z|< \d(z)\}} \frac{h(y) \d^{s_2}(x)}{|x-y|^{N-s_2}}dy,\\
            \\&\displaystyle\le C\int_{\{z\in\O\;;\;\d(x)\le |x-z|< \d(z)\}} \frac{h(y)}{|x-y|^{N-2s_2+s_1}}dy,\\
            \\&\displaystyle\le  C \int_{\{z\in\O\;;\;\d(x)\le |x-z|< \d(z)\}} \frac{h(y)\d^{a(1-s_1)}(y)}{|x-y|^{N-2s_2+s_1} \d^{a(1-s_1)}(y)}dy,\\
            \\&\displaystyle\le C \int_{\O}
            \frac{h(y)\d^{a(1-s_1)}(y)}{|x-y|^{N-(2s_2-s_1-a(1-s_1))}}dy.
        \end{array}\right.
    \end{equation}
    Thus, using again Lemma~\ref{eq:Stein}, it follows that
    $$
    \|J_{3}\|_{L^{\g}(\O)}\le C\|h \d^{a(1-s_1)}\|_{L^{m}(\O)} \qquad \hbox{
        for any } \g\leq \widehat{m} _{s_1,s_2,a}.
    $$
    $\bullet$ Now, let us deal with the last term ~$J_{4}$. From \eqref{eq:First_GPel2} and  since
    $a<\dfrac{s_2}{1-s_1}$ with $s_2>s_1$, we obtain
    \begin{equation}\label{eq:J32(x)el2}
        \left\{\begin{array}{ll}
            J_{4}(x) &\le \displaystyle C \int_{\{z\in\O\;;\;|x-z|\ge \max\{\d(z),\d(x)\}\}} \frac{h(y)\d^{s_2}(y)\d^{s_2}(x)}{\d^{s_1}(x)|x-y|^N}dy,\\
            \\&\le \displaystyle C \int_{\{z\in\O\;;\;|x-z|\ge \max\{\d(z),\d(x)\}\}} \frac{h(y) \d^{a(1-s_1)}(y)\d^{s_2-a(1-s_1)}(y)\d^{s_2-s_1}(x)}{|x-y|^N}dy,\\
            \\&\le \displaystyle C \int_{\{z\in\O\;;\;|x-z|\ge \max\{\d(z),\d(x)\}\}}
            \frac{h(y)\d^{a(1-s_1)}(y)}{|x-y|^{N-(2s_2-s_1-a(1-s_1))}}dy.
        \end{array}\right.
    \end{equation}
    Hence,  we conclude as in the previous estimates. Finally, we have proven that
    \begin{equation*}
        \left\|\frac{w}{\d^{s_1}}\right\|_{L^{\g}(\O)} \leq
        C\|h\d^{a(1-s_1)}\|_{L^{m}(\O)}\qquad \hbox{ for any } \g\leq
        \widehat{m} _{s_1,s_2,a}.
    \end{equation*}
\end{proof}

Next, we state and prove a weighted regularity result for the gradient of the solution to Problem~\eqref{eq:poissonel2}. 


\begin{thm}\label{eq:main_regularity_resultel2}
   Let  $s_1, s_2\in (\frac12,1)$ such that $s_1<s_2$.
     Assume all the hypotheses of Theorem~\ref{eq:regularity_usurdeltas}. Let $w$ be  the unique weak solution to Problem~\eqref{eq:poissonel2}.
    Furthermore, let us set $\widetilde{m}_{s_1,s_2,a}:=\min\{\widehat{m} _{s_1,s_2,a},
    \frac{mN}{(N-m(2s_2-1))_+}\}$ where
    $\widehat{m} _{s_1,s_2,a}:=\frac{mN}{(N-m(2s_2-s_1-a(1-s_1)))_+}$. \\
    We have the following assertions.
    \begin{enumerate}

        \item[\em (i)] If $m<\min\{\frac{N}{2s_2-1},\frac{N}{2s_2-s_1-a(1-s_1)}\}$, then $|\n w|\d^{1-s_1} \in L^{\g}(\O)$ for any
        $\gamma<\widetilde{m}_{s_1,s_2,a}$ and
        \begin{equation}\label{eq:comp_resultel2} 
            \left\||\n w|\d^{1-s_1}\right\|_{L^{\g}(\O)}\leq C\bigg(
            \left\|\dfrac{w}{\d^{s_1}}\right\|_{L^{\widehat{m} _{s_1,s_2,a}}(\O )}+
            \|\mathbb{P}(h\d^{a(1-s_1)})\|_{L^{\gamma}(\O )}\bigg),
        \end{equation}
        where $\mathbb{P}: L^m(\O)\to L^\g(\O)$  is the linear operator defined by
        \begin{equation}\label{eq:Piel2}
            \mathbb{P}(g)(x):=
            \left\{
            \begin{array}{lll}
                \dyle \io\frac{g(y)}{|x-y|^{N-(2s_2-s_1-a(1-s_1))}} dy
                &
                \mbox{\em if } & a\ge 1,\\ \\
                \dyle \io\frac{g(y)}{|x-y|^{N-2s_2+1}} dy &\mbox{\em
                    if
                } & 0\le a<1.
            \end{array}
            \right.
        \end{equation}
        Moreover,
        we have \begin{equation}\label{eq:main000el2}  
            \left\||\n w|\d^{1-s_1}\right\|_{L^{\g}(\O)}\leq C\|h \d^{a(1-s_1)}\|_{L^{m}(\O )}\qquad \forall \g<\widetilde{m}_{s_1,s_2
                ,a}.
        \end{equation}

        \item[\em (ii)]
        If $m=\max\left\{\frac{N}{2s_2-1},\frac{N}{2s_2-s_1-a(1-s_1)}\right\}$, then $|\n
        w|\d^{1-s_1}\in L^\g(\O)$ for any $\g<+\infty$

        \begin{equation}\label{eq:comp_result0122el2} 
            \left\||\n w|\d^{1-s_1}\right\|_{L^{\g}(\O)}\leq C\bigg(
            \left\|\dfrac{w}{\d^{s_1}}\right\|_{L^{\g}(\O )}\,\,\,
            + \|\mathbb{P}(h\d^{a(1-s_1)})\|_{L^{\g}(\O )}\bigg).
        \end{equation}
        \item[\em (iii)]
        If $m>\max\{\frac{N}{2s_2-1},\frac{N}{2s_2-s_1-a(1-s_1)}\}$, then $|\n
        w|\d^{1-s_1}\in L^\infty(\O)$ and

        \begin{equation}\label{eq:comp_result01el2}
            \left\||\n w|\d^{1-s_1}\right\|_{L^{\infty}(\O)}\leq C\bigg(
            \left\|\dfrac{w}{\d^{s_1}}\right\|_{L^{\infty}(\O )}\,\,\,
            + \|\mathbb{P}(h\d^{a(1-s_1)})\|_{L^{\infty}(\O )}\bigg).
        \end{equation}
    \end{enumerate}
    Hence, for $\g<\min\{\frac{N}{N-(2s_2-s_1-a(1-s_1))},
    \frac{N}{N-(2s_2-1)}\}$ fixed, the operator
    $\Gamma:L^1(\O,\d^{a(1-s_1)})\to
    \mathbb{W}^{1,\g}_0(\Omega,\d^{\g(1-s_1)})$ defined by $\Gamma(h)=w$,
    where $w$ is the unique weak solution to \eqref{eq:poissonel2}, is
    compact.
\end{thm}

\noindent\begin{proof} We follow closely the argument used in
\cite{AB01} and \cite{AtmaBirDaouLaam}. For the sake of clarity, we will split the proof into two parts:
    weighted regularity and compactness.

\smallbreak

   \noindent   $\bullet$ \textbf{Weighted Regularity}

    Recall that $w$ is given by
    $$
    w(x)=\io \mathcal{G}_{s_2}(x ,y) h(y)\,dy. $$
    Without loss of generalities, we suppose that $h\geq 0$. Then,
    $$
    |\n w(x)|\le \io |\n_x \mathcal{G}_{s_2}(x ,y)| h(y)\,dy\le \io
    \frac{|\n_x
        \mathcal{G}_{s_2}(x ,y)|}{\mathcal{G}_{s_2}(x ,y)}\mathcal{G}_{s_2}(x ,y)
    h(y)dy. $$ Using the estimate \eqref{eq:Fourth_GPel2}, we get

    \begin{equation*}
        \begin{array}{lll}

            |\n w(x)| &\dyle\le  C\int_{\{z\in\O\;;\;|x-z|<\d(x)\}}
            \frac{\mathcal{G}_{s_2}(x ,y)}{|x-y|}h(y)dy+ \frac{C}{\d(x)}
            \int_{\{z\in\O\;;\;|x-z|\ge \d(x)\}}
            \mathcal{G}_{s_2} (x ,y)h(y) dy.
        \end{array}
    \end{equation*}
    Thus,
    \begin{equation}
        \left\{
        \begin{array}{lll}
            |\nabla w(x)|\d^{1-s_1}(x)&\leq \dyle C\d^{1-s_1}(x) \left(\int_{\{z\in\O\;;\;|x-z|< \d(x)\}}
            \dfrac{h(y)}{|x-y|}\mathcal{G}_{s_2}(x ,y)dy\right)
            \vspace{0.2cm}\\
            &+ C\dyle  \d^{1-s_1}(x) \left(\int_{\{z\in\O\;;\;|x-z|\geq \d(x)\}}
            \dfrac{h(y)}{\d(x)}\mathcal{G}_{s_2}(x ,y)dy\right),
            \vspace{0.2cm}\\&=Q_1(x)+Q_2(x).
        \end{array}
        \right.
    \end{equation}
    --- Let us begin by estimating $Q_2$. We have
    \begin{equation}\label{eq:Q2Rel2}
        Q_2(x)\leq C \dfrac{w(x)}{\d^{s_1}(x)}.
    \end{equation}
    As soon as $\g<\widetilde{m}_{s_1,s_2,a}\le \widehat{m}_{s_1,s_2,a}$, we obtain
    \begin{equation}\label{eq:u_q_d-}
        \io Q^\g_2(x) dx\le C\int_{\O} \dfrac{w^\g(x)}{\d^{s_1\g}(x)} dx\leq
        C\|h
        \d^{a(1-s_1)}\|^\g_{L^{m}(\O )}.
    \end{equation}
    --- Now, we  deal with $Q_1$. We have
    \begin{equation}
        \left\{
        \begin{array}{lll}
            Q_1(x)&= \dyle \d^{1-s_1}(x) \left(\int_{\{z\in\O\;;\;|x-z|< \d(x)\}}  \dfrac{h(y)}{|x-y|}\mathcal{G}_{s_2}(x ,y)dy\right),
            \vspace{0.2cm}\\
            &\le \dyle \d^{1-s_1}(x) \left(\int_{\{z\in\O\;;\;\frac 12 \d(x)\le |x-z|< \d(x)\}}  \dfrac{h(y)}{|x-y|}\mathcal{G}_{s_2}(x ,y)dy\right)
            \vspace{0.2cm}\\
            &+\dyle \d^{1-s_1}(x) \left(\int_{\{z\in\O\;;\;|x-z|<\frac 12 \d(x)\}}  \dfrac{h(y)}{|x-y|}\mathcal{G}_{s_2}(x ,y)dy\right),\\
            &= Q_{11}(x)+Q_{12}(x).
        \end{array}
        \right.
    \end{equation}

    \noindent
    {Regarding $Q_{11}$, we have}
    \begin{equation}\label{eq:Q111Rel2}
        \begin{array}{lll}
            Q_{11}(x)&\leq
            C\dyle \d^{-s_1}(x)\int_{\{z\in\O\;;\;\frac 12 \d(x)\le |x-z|< \d(x)\}}  h(y)\mathcal{G}_{s_2}(x ,y)dy
            \le C\dfrac{w(x)}{\d^{s_1}(x)}.
        \end{array}
    \end{equation}
    As above, since $\g<\widetilde{m}_{s_1,s_2,a}\le \widehat{m}_{s_1,s_2,a}$,
    then
    \begin{equation}\label{eq:Q111el2}
        \io Q^\g_{11}(x) dx\le C\int_{\O}
        \dfrac{w^\g(x)}{\d^{s_1\g}(x)} dx\leq C\|h
        \d^{a(1-s_1)}\|^\g_{L^{m}(\O )}.
    \end{equation}

    \noindent Now, we  estimate the term $Q_{12}$. Recall that for any  $x\in\O$ and $y\in
    \{z\in\O\;;\;|x-z|<\frac 12 \d(x)\}$, we have $ \frac 12 \d(x)\le \d(y)\le
    \frac 32 \d(x)$. Hence, from \eqref{eq:First_GPel2}, we find
    \begin{equation}
        Q_{12}(x)\le C \bigg(\dyle \int_{\{z\in\O\;;\;|x-z|\le \frac 12
            \d(x)\}}\dfrac{h(y)\d^{a(1-s_1)}(y)\d^{1-s_1-a(1-s_1)}(x)}{|x-y|^{N-2s_2+1}}
        dy\bigg).
    \end{equation}
    Therefore, according to the value of $a$, we have
    \begin{equation}\label{eq:Q12(x)el2}
        Q_{12}(x)\le
        \left\{
        \begin{array}{lll}
            C \bigg(\dyle \io\frac{h(y)\d^{a(1-s_1)}(y)}{|x-y|^{N-(2s_2-s_1-a(1-s_1))}} dy\bigg)
            &
            \mbox{ if } & a\ge 1,\\ \\
            C\dyle \bigg( \io\frac{h(y)\d^{a(1-s_1)}(y)}{|x-y|^{N-2s_2+1}} dy\bigg)&\mbox{
                if
            } & 0\le a<1.
        \end{array}
        \right.
    \end{equation}
    Taking into consideration the definition of $\mathbb{P}$, it
    follows that
    $$
    \io Q^\g_{12}(x)  dx \le C\io (\mathbb{P}(h\d^{a(1-s_1)})(x))^\gamma
    dx.
    $$
    Combining the estimates \eqref{eq:Q2Rel2} 
    and \eqref{eq:Q12(x)el2}, it holds that
    \begin{equation}\label{eq:maingwwel2}
        |\nabla w(x)|\d^{1-s_1}(x)\le C\dfrac{w(x)}{\d^{s_1}(x)}+
        \mathbb{P}(h\d^{a(1-s_1)})(x).
    \end{equation}
    Hence \eqref{eq:comp_resultel2} follows.

    \

    \noindent Finally, by using Theorem \ref{eq:regularity_usurdeltas} and Lemma \ref{eq:Stein}, we deduce that from \eqref{eq:maingwwel2},

    $$
    \left\||\n w|\d^{1-s_1}\right\|_{L^{\g}(\O)}\leq C\|h \d^{a(1-s_1)}\|_{L^{m}(\O )},
    $$
    for any $\g<\widetilde{m}_{s_1,s_2,a}$.

    \medbreak

    \noindent$\bullet$ \textbf{Compactness\\}
    Now, let us prove that $\Gamma$ is compact. So,
    let $\{h_n\}_n\subset L^1(\Omega, \d^{a(1-s_1)})$ be a
    bounded sequence  and set $w_n=\Gamma(h_n)$. Then, $w_n$ solves
    \begin{equation}
        \left\{
        \begin{array}{rclll}
            (-\Delta )^{s_2} w_n=h_n &\hbox{   in   } & \Omega,\\ w_n=0   &\hbox{
                in }& \mathbb{R}^N\setminus\Omega.
        \end{array}
        \right.
    \end{equation}
    For $q<\frac{N}{N-(2s_2-s_1-a(1-s_1))}$ fixed, by Theorem~\ref{eq:main_regularity_resultel2}, there exists $C>0$ independent of $n$
    such that
    \begin{equation}\label{eq:borneinLql2}
        \int_{\O} \dfrac{|w_n(x)|^q}{\d^{s_1q}(x)} dx\leq C\|h_n
        \d^{a(1-s_1)}\|^q_{L^{1}(\O )}\le C.
    \end{equation}
    Hence, using the estimate~\eqref{eq:comp_resultel2} and choosing $q$ such that  $\g<q<\frac{N}{N-(2s_2-s_1-a(1-s_1))}$, we conclude that
    \begin{equation} \left\{
        \begin{array}{lll}
            \left\||\nabla w_n|\d^{(1-s_1)}\right\|_{L^\g(\O)} &\leq &  C\left(
            \left\|\dfrac{w_n}{\d^{s_1}}\right\|_{L^{q}(\O )}+ \|\mathbb{P}(h_n\d^{a(1-s_1)})\|_{L^\g(\O)}\right),\\ &\le & C\|h_n
            \d^{a(1-s_1)}\|_{L^{1}(\O )}\le C.
        \end{array}
        \right.\end{equation}
    Taking into consideration that  $|\n (w_n\d^{1-s_1})|\le \d^{1-s_1}|\n
    w_n|+ (1-s_1)\dfrac{|w_n(x)|}{\d^{s_1}(x)}$, we deduce that the sequence
    $\{w_n\d^{1-s_1}\}_n$ is bounded in the space $\mathbb{W}^{1,\g}_{0}(\O)$ for any
    $\g<\min\{\frac{N}{N-(2s_2-s_1-a(1-s_1))},\frac{N}{N-(2s_2-1)}\}$. Then, up
    to a subsequence, we get the existence of a measurable function
    $w$ such that $w\d^{1-s_1}\in \mathbb{W}^{1,\g}_0(\Omega)$,
    $w_n\d^{1-s_1}\rightharpoonup w\d^{1-s_1}$ weakly in
    $\mathbb{W}^{1,\g}_{0}(\Omega)$ and $w_n\to w$ \textit{a.e.} in $\O$.
    \noindent Thanks to \eqref{eq:borneinLql2},  $\left\{\dfrac{w_n}{\d^{s_1}}\right\}_n$ is bounded in $L^q(\O)$ for
    any $q<\frac{N}{N-(2s_2-s_1-a(1-s_1))}$. Therefore,  we
    conclude that $\dfrac{w}{\d^{s_1}}\in L^q(\O)$ for any
    $q<\frac{N}{N-(2s_2-s_1-a(1-s_1))}$ by Fatou's Lemma. Finally, since $\O$ is a
    bounded {open subset of $\mathbb{R}^N$}, by using Vitali's Lemma we obtain that
    $$
    \dfrac{w_n}{\d^{s_1}}\longrightarrow \dfrac{w}{\d^{s_1}}\mbox{  strongly
        in }L^{q}(\O),$$ for any $q<\frac{N}{N-(2s_2-s_1-a(1-s_1))}$.
    \noindent Now, setting $v_{ij}=w_i-w_j$ and  $h_{ij}=h_i-h_j$, it holds that $v_{ij},
    \dfrac{v_{ij}}{\d^{s_1}}\longrightarrow 0$ \textit{a.e.} in $\O$ as
    $i,j\longrightarrow +\infty$ and
    $\dfrac{v_{ij}}{\d^{s_1}}\longrightarrow 0$ strongly in $L^{q}(\O)$ for
    any $q<\frac{N}{N-(2s_2-s_1-a(1-s_1))}$.
    \noindent Since $v_{ij}$ solves  {the Problem}
    \begin{equation}
        \begin{cases}
            (-\Delta )^{s_2} v_{ij}=h_{ij} &\hbox{   in   }  \Omega,\\ v_{ij}=0
            &\hbox{   in   } \mathbb{R}^N\setminus\Omega.
        \end{cases}
    \end{equation}
    From \eqref{eq:comp_resultel2}, we get
    $$
    \begin{array}{lll}
        \left\||\nabla v_{ij}|\d^{1-s_1}\right\|_{L^\g(\O)} &\leq &  C\left(
        \left\|\dfrac{v_{ij}}{\d^{s_1}}\right\|_{L^{q}(\O )}+ \|\mathbb{P}(h_{ij}\d^{a(1-s_1)})\|_{L^\g(\O)}\right).
    \end{array}
    $$
    Recall that $\dfrac{v_{ij}}{\d^{s_1}}\to 0\mbox{  as  }i,j\to +
    \infty$ strongly in $L^q(\O)$ for any $q<\frac{N}{N-(2s_2-s_1-a(1-s_1))}$. Therefore, we deduce from \cite[Theorem 8.1]{KZP}
    that the operator
    $\mathbb{P}: L^1(\O)\to L^\g(\O)$ is compact for any
    $\g<\min\{\frac{N}{N-(2s_2-s_1-a(1-s_1))},
    \frac{N}{N-(2s_2-1)}\}$.
    Hence, up to an other subsequence, we conclude that
    $w_n\longrightarrow w$ strongly in $\mathbb{W}^{1,\g}_0(\Omega,\d^{\g(1-s_1)})$. This ends the proof.
\end{proof}

\section{Existence and nonexistence results}\label{eq:Existenc_Nonexistence_S}
As mentioned in the introduction, our main objective is to {investigate} the existence and nonexistence of solutions to
   	$$ (S)\hspace*{.25cm}
    \left\{
    \begin{array}{rclll}
        (-\Delta)^{s_1} u &= & |\nabla v|^{q}+\l f & \text{ in }&\Omega , \\
        (-\Delta)^{s_2} v &= & |\nabla u|^{p}+\mu g & \text{ in }&\Omega , \\
        u =v&=&0 &\hbox{  in }& \mathbb{R}^N\setminus\Omega,
    \end{array}%
    \right.
    $$
 where  $\frac 12<s_1< s_2<1$. \textbf{\em We emphasize that this assumption will hold throughout the rest of this paper}, and that the case $s_2< s_1$ could be handled in the same manner.
    
    \smallbreak
\noindent In the first subsection, we will establish the existence of a solution to System $(S)$ under natural assumptions on the data. Then, we will establish nonexistence results in the second
subsection.
\smallbreak
\noindent To begin with, let us specify what we mean by a (weak) solution to System $(S)$.
\begin{defi}\label{eq:def:faible}
	Let $(f, g)\in (L^1(\O)^+)^2$. We say that $(u,v)\in \mathbb{W}^{1,p}_0(\O)\times \mathbb{W}^{1,q}_0(\O)$    is a (weak) solution to System $(S)$ if 
	for any $(\varphi, \psi)\in \mathbb{X}_{s_1}(\Omega)\times \mathbb{X}_{s_2}(\Omega)$ (see Definition \ref{eq:def1o}), we have
	\begin{equation}\label{eq:eq1fel2}
		\left\{\begin{array}{ll}
			\dyle\io u (-\Delta)^{s_1}\varphi\, dx
			&=\dyle\io |\n v|^q\varphi\, dx+\l\dyle\io f \varphi\, dx \\ \dyle\io v(-\Delta)^{s_2}\psi \, dx&=\dyle\io| \nabla u|^p\psi\, dx+\mu\dyle\io g\psi\, dx.
		\end{array}
		\right.
	\end{equation}
\end{defi}

\subsection{Existence results}

In this subsection, we will state and prove our  existence theorems in the following two cases:\vspace*{0.2cm}
\\
1) $p,q\geq 1$ and $pq>1$;\vspace*{0.2cm}\\
2)  $p=q=1$.

\smallbreak

\noindent As we shall see, the technical assumptions underlying the regularity results of Section~\ref{eq:sec:regu-comp} naturally shape the hypotheses of our main theorems.
For the ease of exposition, we introduce the following notations, which will be extensively used in the formulation of said hypotheses.
	\begin{itemize}
	\item Let $\kappa$ denote either $p$ or $q$;
	\item $c_\kappa:=\dfrac{\kappa+1}{\kappa+2}$ and $\widehat{\kappa}_{i,s}:=is-\kappa(1-s)$ for $i\in\{1,2\}$;
		\item $\widehat{\kappa}_{s_1,s_2}:=2s_2-s_1-\kappa(1-s_1)$;
	\item $\widehat{m}_{s_1,\kappa}:=\left\{ \begin{array}{llll}
	\displaystyle\frac{mN}{N-m\widehat{\kappa}_{1,s_1}}& \text{if } \kappa >1,
	\\ 
	\displaystyle\frac{mN}{N-m(2s_1-1)} & \text{if } \kappa =1,	
	\end{array}\right.$ and $\overline{m}_{s_1,s_2,\kappa}:= \displaystyle\frac{\widehat{m}_{s_1,\kappa}}{1+ \widehat{m}_{s_1,\kappa}(s_2-s_1)}$;
	\item $\widehat{\s}_{s_2}:=\displaystyle\frac{\s N}{N-\s(2s_2-1)}$  and $ \widehat{\s}_{s_1,s_2,\kappa}:=\displaystyle\frac{\s N}{N-\s\widehat{\kappa}_{s_1,s_2}}$.
\end{itemize}

\subsubsection{Case:   $p,q\geq 1$ and $pq>1$}
	\noindent In light of what we will see in the proof, the assumptions  {$1\leq q<\frac{s_1}{1-s_1}$ and $1\leq p<\frac{s_2}{1-s_1}$}
are necessary to establish existence and apply our regularity results. These requirements  naturally feature in the statement of our first main theorem presented next.

\begin{thm}\label{eq:existence_theorem_for_Sel2}
   Assume that  {$1\leq q<\dfrac{s_1}{1-s_1}$ and $1\leq p<\dfrac{s_2}{1-s_1}$}  such that $pq>1$. Suppose that $(f,g) \in L^{m}(\O)^+\times L^{\s}(\Omega)^+$ where $(m,\s)\in [1,+\infty)^2$ satisfies one of the following set of conditions:
    \begin{equation}
        \left\{
        \begin{array}{llll}\label{eq:condi02022}
          c_q<s_1<1,\,\, c_p <s_2<1, \vspace{0.2cm}\\
            1<m<\dfrac{N}{\widehat{q}_{1,s_1}},\,\,\,\, 1<\s<\dfrac{N}{\widehat{p}_{s_1,s_2}},\\
            p\s<\overline{m}_{s_1,s_2,q}\quad \text{and} \quad
            qm < \widehat{\s}_{s_1,s_2,p},\,
        \end{array}%
        \right.
    \end{equation}
    or
    \begin{equation}
        \left\{
        \begin{array}{llll}\label{eq:condi0012022}
            1/2<s_1<c_q,\,\,\, 1/2<s_2<c_p,

            \\
            \dfrac{1}{\widehat{q}_{2,s_1}}<m<\dfrac{N}{\widehat{q}_{1,s_1}},\,\,\,\dfrac{1}{\widehat{p}_{2,s_2}}<\s<\dfrac{N}{\widehat{p}_{s_1,s_2}} ,\vspace{0.2cm}\\
            p\s<\overline{m}_{s_1,s_2,q}\quad \text{and} \quad qm<\widehat{\s}_{s_1,s_2,p},
        \end{array}%
        \right.
    \end{equation}
    or
    \begin{equation}
        \left\{
        \begin{array}{llll}\label{eq:condi0002022}
           c_q<s_1<1,\,\, 1/2<s_2<c_p, \,\,q<p,\\
            1<m<\dfrac{N}{\widehat{q}_{1,s_1}},\,\, \dfrac{1}{\widehat{p}_{2,s_2}}<\s<\dfrac{N}{\widehat{p}_{s_1,s_2}}, \vspace{0.2cm}\\
            p\s<\overline{m}_{s_1,s_2,q}\quad \text{and} \quad
            qm < \widehat{\s}_{s_1,s_2,p},
        \end{array}%
        \right.
    \end{equation}
    or
    \begin{equation}
        \left\{
        \begin{array}{llll}\label{eq:condi00i2022}
           1/2<s_1<c_q, c_p<s_2<1,\\\dfrac{1}{\widehat{q}_{2,s_1}}<m<\dfrac{N}{\widehat{q}_{1,s_1}},\,\,\, 1<\s<\dfrac{N}{\widehat{p}_{s_1,s_2}} , \vspace{0.2cm}\\p\s<\overline{m}_{s_1,s_2,q}\quad \text{and} \quad
            qm < \widehat{\s}_{s_1,s_2,p},
        \end{array}%
        \right.
    \end{equation}
    or

    \begin{equation}\label{eq:condi0112022}
        m\ge \dfrac{N}{\widehat{q}_{1,s_1}},\hspace{3mm}\,\s>\dfrac{qmN}{N+qm\widehat{p}_{s_1,s_2}} \text{ and  }\,\,\, p\s<\dfrac{1}{s_2-s_1},
    \end{equation}
    or
   {\begin{equation}\label{eq:condi022022}
            \s \ge \dfrac{N}{\widehat{p}_{s_1,s_2}}\; \text{ and  } \; m>\dfrac{p\s N}{N+p\s(\widehat{q}_{1,s_1}-N(s_2-s_1))}.\quad
    \end{equation}}

  Then, there exists $A>0$ such that for any $(\lambda, \mu)\in \Pi$
    where
    \begin{equation}\label{eq:pi1000-}
        \Pi : =\left\{(\lambda, \mu) \in (0,+\infty)\times (0,+\infty) \;\;;\,\, \lambda^p\|f\|_{L^m(\O)}^p+\mu \|g\|_{L^\s(\O)} \leq A\right\},
    \end{equation}
    System $(S)$ has a nonnegative solution $(u,v)$.
    Moreover, {$(u,v)\in
        \mathbb{W}^{1,\theta_1}_0(\O,\d^{\theta_1(1-s_1)})\times \mathbb{W}^{1,
            \theta_2}_0(\O,\d^{\theta_2(1-s_1)})$ for any
        $\theta_1<\overline{m}_{s_1,s_2,q}$ and
        $\theta_2<\widehat{\s}_{s_1,s_2,p}$}.
 \end{thm} \noindent 

\noindent\textbf{\noindent  Comments.}  Before giving the proof, let us discuss the assumptions of this theorem.      	        	
  %
  	\smallbreak
\noindent {\it First,} we may take $s_2=1$ under the stated assumptions.
 This enables us to handle both local and nonlocal operators within the same system.
   \smallbreak
  \noindent  {\it Second,} in Assumption \eqref{eq:condi0002022}, the complementary case $p<q$ could be treated analogously.  
        \smallbreak
\noindent  {\it Third,} 
to illustrate the assumptions \eqref{eq:condi02022}, \eqref{eq:condi0112022} and \eqref{eq:condi022022}, let us take a concrete case, namely $s_1=\frac34$ and $s_2=\frac9{10}$. 
 \begin{itemize}
  \item[(1)]  	Initially, let us make explicit the size conditions on $(p, q)$ for a given $(m,\sigma)$.
         \begin{itemize}  	
    {\item[(i)] If $m=\sigma=N=2$, then  \eqref{eq:condi02022} is  fulfilled for $q\in [1,\frac95)$ and $p\in [1,\frac{20}{5q+11})$ with $pq>1$.}
    {\item[(ii)] $m=3$ and $\sigma=N=2$, then  \eqref{eq:condi02022} is  satisfied for $q\in [1,2)$ and $p\in [1,\frac{60}{15q+13})$ with $pq>1$.}
     \item[(iii)] $m=\sigma=N=3$, then  \eqref{eq:condi0112022} is  met for $q\in (3,\frac{18}{5})$ and $p\in [1,\frac{q+20}{5q})$.
     \item[(iv)] $m=\sigma=+\infty$ and $N\geq 2$, then  \eqref{eq:condi022022} holds for $q\in (1,3]$ and $p=1$.
         \end{itemize}
\item[(2)] Reciprocally,
        \begin{itemize}  	
      \item[(i)] If $p=\frac32$, $q=1$ and $N=2$, then  \eqref{eq:condi02022} is  satisfied if $
      \begin{cases}
      	\sigma \in(1,\frac{80}{47}], \;\;  m \in(\frac{30\sigma}{3\sigma+20},\frac{80\sigma}{80-27\sigma})
      	\\
      	\text{or} \;\;\sigma \in(\frac{80}{47},\frac{80}{27}),\;\;   m\in (\frac{30\sigma}{3\sigma+20},4).
      \end{cases}
      $
     \item[(ii)]  If $p=2$, $q=\frac32$ and $N=2$, then  \eqref{eq:condi02022} is  fulfilled if  $
     \begin{cases}
     	\sigma \in(\frac{20}{9},\frac{5}{2}], \;\;  m \in(\frac{80\sigma}{3\sigma+40},\frac{80\sigma}{3(40-11\sigma)})
     	\\
     	\text{or} \;\;\sigma \in(\frac{5}{2},\frac{10}{3}),\;\;  m \in(\frac{80\sigma}{3\sigma+40},\frac{16}{3}).
     \end{cases}
     $
       \item[(iii)]  If $p=1$, $q=\frac32$ and $N=3$, then  \eqref{eq:condi02022} holds if $
       \begin{cases}
       	\sigma \in(\frac{12}{9},\frac{20}{7}], \;\;  m \in(\frac{40\sigma}{40-\sigma},\frac{10\sigma}{15-4\sigma})
       	\\
       	\text{or} \;\;\sigma \in(\frac{20}{7},\frac{15}{4}),\;\;  m \in(\frac{40\sigma}{40-\sigma},8).
       \end{cases}
       $
       \end{itemize}
\end{itemize}
    \smallbreak
    \noindent{\it Finally,} within the assumptions of the theorem, setting  $s_1=s_2$ allows one to recover the result established in~\cite{AtmaBirDaouLaam}.

\medbreak

\noindent{\bf Proof of Theorem \ref{eq:existence_theorem_for_Sel2}.} We
follow closely the arguments used in \cite{AB01, AtmaBirDaouLaam,
    AttaBentLaam}. Also, we will give the proof under the assumptions
\eqref{eq:condi02022} and \eqref{eq:condi0012022}. The other cases
follow by using the same reasoning. Since we will use  Schauder's
fixed point Theorem, and in order to make our demonstration  clear
for the readers, we will proceed by three steps.
\\
\textbf{First case:} assume that the assumption \eqref{eq:condi02022} is fulfilled.
\\
Let us define the  following function
$$
\Upsilon(\a):=\a^{\frac{1}{pq}} -\widetilde{C}\a\qquad \forall \a>0,
$$
where $\widetilde{C}$ is a  positive constant depending only on the data and which we will specify later.\\
Since $pq>1$, then there exists $\a_0>0$ such that $\Upsilon(\a_0) =0$, $\Upsilon(\a) >0$   for any $\a\in (0,\a_0)$ and $\Upsilon(\a) <0$  for any $ \a\in (\a_0, +\infty)$. Hence, there exist two positive constants $\ell>0$ and $\L^*>0$  such that
$$
\max_{\a> 0}\Upsilon(\a)=\Upsilon(\ell)=\L^*.
$$
Thus,
$$
\ell^{\frac{1}{pq}}=\widetilde{C}\left(\ell+\frac{\L^*}{\widetilde{C}}\right).
$$
Let $\ell>0$ be fixed as above. Since $(\l,\mu)\in \Pi$ with the choice $A=\frac{\L^*}{\widetilde{C}}$, we have
\begin{equation}\label{eq:relation-entre-l-f-g-lambda, alpha}
    \widetilde{C}\bigg(\ell+\l^p||f||_{L^{m}(\O)}^p+\mu||g||_{L^{\s}(\O)}\bigg)\le \ell^{\frac{1}{pq}}.
\end{equation}
\noindent Since $qm<\widehat{\s}_{s_1,s_2,p}=\frac{\s N}{N-\s(2s_2-s_1-p(1-s_1))}$, then there exists $r>0$ such that $qm<r<\widehat{\s}_{s_1,s_2,p}$. Furthermore, let us consider the set ${H}$ defined by

\begin{equation}\label{eq:sett2}
    {H}:=\Big\{\varphi\in \mathbb{W}^{1,1}_0(\O) \,\, ; \,\, \varphi\, \d^{1-s_1}\in \mathbb{W}^{1,r}_0(\O)\mbox{  and  } \bigg(\io |\n (\varphi\, \d^{1-s_1})|^{r}  dx\bigg)^{\frac{1}{r}}\le \ell^{\frac{1}{pq}}\Big\}.
\end{equation}

\noindent   It is clear  that ${H}$ is a convex and closed bounded subset of $\mathbb{W}^{1,1}_0(\O)$.\\

Let $\varphi\in {H}$ and let us set
$\widehat{\varphi}:=\varphi\, \d^{1-s_1}$.
Then, by using Hardy's inequality (Theorem \ref{eq:HardyInequalityel2}) with $\gamma=0$, we get
$$
\io \frac{|\widehat{\varphi}|^{r}}{\d^{r}} dx \le C(\O)\io |\n \widehat{\varphi}|^{r}  dx.
$$
Thus,
$$
\io \frac{|\varphi|^{r}}{\d^{rs_1}} dx \le C(\O)\io |\n \widehat{\varphi}|^{r}  dx.
$$
Since $|\n \d(x)|=1$ \textit{a.e.} in $\O$, we deduce that
$$
|\n \varphi|^{r}\d^{r(1-s_1)}=|\n (\varphi  \d^{1-s_1})- \varphi\n
\d^{1-s_1}|^{r}\le C_1 |\n \widehat{\varphi}|^{r}+
C_2\frac{|\varphi|^{r}}{\d^{rs_1}}.
$$
{Therefore, using Hardy's inequality}, we obtain $|\n
\varphi|^{r}\d^{r(1-s_1)}\in L^1(\O)$ and
$$
\bigg(\io |\n \varphi|^{r}\, \d^{r(1-s_1)} dx\bigg)^{\frac{1}{r}}\le C \ell^{\frac{1}{pq}}.
$$
We observe that, if  $\varphi \in {H}$,  then $|\n \varphi|^{a} \d^{a(1-s_1)}\in L^1(\O) $ for any $a\le r$.  In particular,  $|\n \varphi|^{q} \d^{q(1-s_1)}\in L^m(\O) $.\\

\noindent {\bf Step 1.} We claim that if $\varphi\in {H}$, then there exists $\b<2s_1-1$ such that $|\n \varphi|^{q} \d^{\b}\in L^1(\O)$.
Indeed, let $\varphi\in {H}$, then by  using  H\"older's inequality,  we get
\begin{equation}  \left\{
    \begin{array}{llll}
        \dyle\io|\n \varphi|^{q}\d^\b  dx&=&\dyle\io|\n \varphi|^{q}\d^{q(1-s_1)}\d^{\beta-q(1-s_1)} dx,\\
        &\dyle\le &\dyle\textbf{}\bigg(\io|\n \varphi|^{qm}\d^{q(1-s_1)m} dx\bigg)^{\frac{1}{m}}
        \bigg(\io \d^{(\beta-q(1-s_1))m'} dx\bigg)^{\frac{1}{m'}},\\
        &\dyle\le&\dyle C \bigg(\io \d^{(\beta-q(1-s_1))m'} dx\bigg)^{\frac{1}{m'}}.
    \end{array}
    \right.
\end{equation}
Since $ \dfrac{q+1}{q+2}<s_1$ (from the set of conditions  \eqref{eq:condi02022}), then there exists $\beta>0$ such that $q(1-s_1)<\beta<2s_1-1$.
Therefore,
$|\nabla \varphi|^{q}\d^\b+\l f \in L^1(\O)$. By Proposition \ref{eq:Existence_with_weightel2}, we get the existence of a  unique weak  solution  $u$ to {the following problem}

\begin{equation}\label{eq:P01el2}
    \left\{
    \begin{array}{rclll}
        (-\Delta)^{s_1} u &= & |\nabla \varphi|^{q}+\l f & \text{ in }&\Omega , \\
        u &=& 0 &\hbox{  in } &\mathbb{R}^N\setminus\Omega.
    \end{array}%
    \right.
\end{equation}%

\noindent   Moreover,  $ u\in \mathbb{W}^{1,\a}_0(\O)$ for any $\a<\dfrac{N}{N-2s_1+1+\beta}$.\\

Now, we claim that $|\n u|^{p} \d^{\b_1}\in L^1(\O)$ for any $\beta_1<2s_2-1$. Indeed,
by Theorem \ref{eq:regubeta1} and by taking  $s=s_1$, $a=q$ and $h=|\n \varphi|^q+\l f$, we have

$$\left\||\n u|\d^{1-s_1}\right\|_{L^{\g}(\O)}\leq C\bigg(\left\||\n \varphi|\d^{1-s_1}\right\|^q_{L^{qm}}+\l \|f\|_{L^m(\O)}\bigg),\qquad  \hbox{      for any         }\,\, \g\leq \widehat{m}_{s_1,q}.$$
Since   $qm<r$, then by H\"older's inequality we get

\begin{equation}\label{eq:uuel2}
    \left\{
    \begin{array}{llll}
        \left\||\n u|\d^{1-s_1}\right\|_{L^{\g}(\O)}&\leq& C\bigg(\left\||\n \varphi|\d^{1-s_1}\right\||^q_{L^{r}(\O)}+\l \|f\|_{L^m(\O)}\bigg),\vspace{0.2cm}\\
        &\leq& C\bigg(\ell^{\frac{1}{p}}+\l \|f\|_{L^m(\O)}\bigg),
    \end{array}
    \right.
\end{equation}
for any $ \g\leq \widehat{m}_{s_1,q}$. Thus,  $|\n u|^p \d^{p(1-s_1)}\in L^{\s}(\O)\subset L^{1}(\O)$.  Since $p\s<\overline{m}_{s_1,s_2,q}$, then by using H\"older's inequality, we  obtain $|\n u|^p \d^{p(1-s_2)}\in L^{\s}(\O)$.\\

\noindent   On the other hand, as   $\dfrac{p+1}{p+2}<s_2$,  we get by using the same computations as    above,   the existence of  $\b_1<2s_2-1$, such that $|\n u|^p \d^{\b_1}\in L^1(\O)$.

\noindent   Thus, using again Proposition \ref{eq:Existence_with_weightel2}, we get the existence of a weak solution $v$  to the following {problem}

\begin{equation}\label{eq:P02el2}
    \left\{
    \begin{array}{rclll}
        (-\Delta)^{s_2} v &= & |\nabla u|^{p}+\mu g & \text{ in }&\Omega , \\
        v &=& 0 &\hbox{  in }& \mathbb{R}^N\setminus\Omega,\\
    \end{array}%
    \right.
\end{equation}%
such that  $v\in \mathbb{W}^{1,\a_1}_0(\O)$ for any $\a_1<\dfrac{N}{N-2s_2+1+\b_1}$.

\noindent Therefore,   the operator
$$T:  {H}\longrightarrow \mathbb{W}^{1,1}_0(\O)\  \  $$
$$\varphi \longmapsto T(\varphi)=v$$
is well defined. Moreover, if $v$ is a fixed point of $T$, then $(u,v)$ is a weak solution to  System $(S)$.

\vspace{0.2cm}
\noindent {\bf Step 2.} In this step, we want to show that $T(H)\subset H$.
\\
Indeed, by taking  $h= |\n u|^p+\mu g$ and  by using Theorem \ref{eq:main_regularity_resultel2} {with $a=p$}, it holds that
\begin{equation}\label{eq:esti1o}
    \left\||\n v|\d^{1-s_1}\right\|_{L^{\nu}{(\O)}}\leq C\bigg(\left\||\n u|\d^{1-s_1}\right\|^{p}_{L^{p\s}(\O)}+{\mu}\|g\|_{L^{\s}(\O)}\bigg),\qquad \hbox{      for any         }\quad \nu\leq \widehat{\s}_{s_1,s_2,p}.
\end{equation}
Since, $p\s<\overline{m}_{s_1,s_2,q}<\widehat{m}_{s_1,q}$ and  thanks to estimate  \eqref{eq:uuel2}, it follows that
\begin{equation}\label{eq:alphaa}
    \left\| |\n v| \d^{1-s_1}\right\|_{L^{\nu}(\O)}\le C\bigg(\left\|| \n \varphi|\d^{1-s_1}\right\|^{pq}_{L^{r}(\O)}+\mu \|g\|_{L^\s(\O)}+ \l^p \|f\|^p_{L^m(\O)}\bigg).
\end{equation}
Since $r<\widehat{\s}_{s_1,s_2,p}$, then by choosing $\nu=r$ in the previous inequality, we get
\begin{equation*}
    \left\| |\n v| \d^{1-s_1}\right\|_{L^{r}(\O)}\le C\bigg(\ \left\||\n \varphi|\d^{1-s_1}\right\|^{pq}_{L^{r}(\O)}+\mu \|g\|_{L^\s(\O)}+ \l^p\|f\|^p_{L^m(\O)}\bigg).
\end{equation*}
Recall that $\varphi\in H$, thus
$$
\left\| |\n v| \d^{1-s_1}\right\|_{L^{r}(\O)}\le C\bigg(\ell+\mu \|g\|_{L^{\s}(\O)}+\l^p \|f\|^p_{L^{m}(\O)}\bigg).
$$
By choosing $\widetilde{C}=C$ and taking into consideration the definition of $\ell$,
we conclude that\\
$$\left\| |\n v| \d^{1-s_1}\right\|_{L^{r}(\O)} \leq \ell^{\frac{1}{pq}}.$$Thus,  $v\in H$ and therefore $T(H)\subset H$.\\

\noindent {\bf Step 3.} Continuity and Compactness of $T$.
\vspace{0.2cm}\\
$\bullet$  \textbf{\underline{Continuity} }\vspace{0.2cm}\\
-- Let $\{\varphi_n\}_n\subset H$ and $\varphi\in H$
be such that
$\varphi_n \to \varphi \,\,\text{ strongly in } \mathbb{W}^{1,1}_0(\O)$. Moreover, let us define $v_n=T(\varphi_n)$ and $v=T(\varphi)$. Then, $(u_n,v_n)$  and $(u,v)$ satisfy
\begin{equation}\label{eq:PP01el2}
    \left\{
    \begin{array}{cclll}
        (-\Delta)^{s_1} u_n &= & |\n\varphi_{n}|^q +\l f & \text{ in }&\Omega, \\
        (-\Delta)^{s_1} u &= & |\n\varphi|^q+\l f & \text{ in }&\Omega, \\
        u_n= u &=& 0  & \hbox{  in } &\mathbb{R}^N\setminus\Omega,\\
    \end{array}
    \right.
\end{equation}
and
\begin{equation}\label{eq:PP02el2}
    \left\{
    \begin{array}{cclll}
        (-\Delta)^{s_2} v_n & = & |\nabla u_n|^{p}+\mu g
        & \text{ in }&\Omega,\\
        (-\Delta)^{s_2} v & = & |\nabla u|^{p}+\mu g
        & \text{ in }&\Omega,\\
        v_n= v&=&0  & \hbox{  in }& \mathbb{R}^N\setminus\Omega.
    \end{array}
    \right.
\end{equation}

\noindent   Thus, by using Proposition \ref{eq:Existence_with_weightel2}, we get
$$\left\|\n u_n-\n u\right\|_{L^\a(\O)}\le C\left\| |\n\varphi_n|^q-|\n \varphi|^q\right\|_{L^1(\O,\d^{\b})},
$$
for any $\a<\frac{N}{N-2s_1+1+\b}$ where $\b<2s_1-1$.\\

\noindent
Since  $\varphi_n \to \varphi \text{ in } \mathbb{W}^{1,1}_0(\O)$  and $\{\varphi_n\}_n\subset H$, then by Fatou's Lemma, we get $$\bigg(\io |\n \varphi|^{r}\, \d^{r(1-s_1)} dx\bigg)^{\frac{1}{r}}\le {C} \ell^{\frac{1}{pq}}.$$
\noindent
Hence, for $a\in (1,r)$   fixed and due to Proposition \ref{eq:interpolation}, we get

\begin{equation}
    \left\{
    \begin{array}{rclll}
        \bigg(\dyle\io|\n \varphi_n-\n \varphi|^{a}\d^{a(1-s_1)} dx\bigg)^{\frac{1}{a}} &\le & \bigg(\dyle\io|\n \varphi_n-\n \varphi|\d^{1-s_1} dx\bigg)^{\eta}\\ &\times &\bigg(\dyle\io|\n \varphi_n-\n \varphi|^{r}\d^{r(1-s_1)} dx\bigg)^{\frac{1-\eta}{r}}, \\
        & \le & C \bigg(\dyle\io|\n \varphi_n-\n \varphi|\d^{1-s_1} dx \bigg)^{\eta}\to 0\mbox{  as   }n\to +\infty,
    \end{array}
    \right.\end{equation}
where $\eta=\dfrac{r-a}{a(r-1)}$.\\
Since $qm<r$, thus by choosing $a=qm$ we can find, as in Step 1,  $\b<2s_1-1$ such that
$$
\begin{array}{lll}
    \dyle\io|\n \varphi_n-\n \varphi|^{q}\d^{\b} dx &\le & C\bigg(\dyle\io|\n \varphi_n-\n \varphi|^{qm}\d^{qm(1-s_1)} dx\bigg)^{\frac{1}{m}}\to 0\mbox{  as   }n\to +\infty{.}
\end{array}
$$

\noindent
Thus, it follows that
$$
\left\|\n u_n-\n u\right\|_{L^\a(\O)}\to 0\mbox{  as  }n\to +\infty,
$$
for any $\a<\dfrac{N}{N-2s_1+1+\b}$. In particular,  $u_n\to u$ strongly in $\mathbb{W}^{1,1}_0(\O)$.
\vskip1mm
\noindent Now  going back to \eqref{eq:uuel2}, we deduce that
\begin{equation}\label{eq:RR10-}
    \left\| |\n u_n| \d^{1-s_1}\right\|_{L^{\widehat{m}_{s_1,q}}(\O)}\le C\bigg(\ell^{\frac 1p}+\l \|f\|_{L^{m}(\O)}\bigg) .
\end{equation}
Thus, $\{|\n u_n|\d^{1-s_1}\}_n$ is bounded in $L^{\widehat{m}_{s_1,q}}(\O)$. \\
\noindent
Now, we claim that
\begin{equation}\label{eq:boundel2}
    \{|\n u_n|\d^{1-s_2}\}_n  \hbox {      is bounded in      } L^{\g}(\O) \hbox  {      for any      }\g<\overline{m}_{s_1,s_2,q}.\\
\end{equation}
\noindent
Indeed, let $\g<\overline{m}_{s_1,s_2,q}$, then by using the fact  that $\{|\n u_n|\d^{1-s_1}\}_n$ is bounded in $L^{\widehat{m}_{s_1,q}}(\O)$ and {by applying} H\"older's inequality,  we get
$$
\begin{array}{lll}
    \dyle\io|\n u_n|^{\g}\d^{\g(1-s_2)} dx &\le & \bigg(\dyle\io|\n u_n|^{\widehat{m}_{s_1,q}}\d^{\widehat{m}_{s_1,q}(1-s_1)} dx\bigg)^{\frac{\g}{\widehat{m}_{s_1,q}}}\bigg(\dyle\io\frac{ dx}{\d^{\frac{\widehat{m}_{s_1,q}\g(s_2-s_1)}{\widehat{m}_{s_1,q}-\g}}}\bigg)^{1-\frac{\g}{\widehat{m}_{s_1,q}}}<+\infty.
\end{array}
$$
Here, we have used the fact that  $\g<\overline{m}_{s_1,s_2,q}$. Then, the claim follows.
\noindent
Therefore, by Fatou's Lemma, it follows that $$\io |\n u|^{\g}\, \d^{\g(1-s_2)} dx<+\iy.$$
Let $\g<\overline{m}_{s_1,s_2,q}$ be fixed such that $p\s<\g$, then by using  Proposition \ref{eq:interpolation} we get
\begin{equation}\label{eq:Strong_conv_uel2}
    \left\{     \begin{array}{lll}
        \bigg(\dyle\io|\n u _n-\n u|^{p\s}\d^{p\s (1-s_2)} dx\bigg)^{\frac{1}{p\s}} &\le & \bigg(\dyle\io|\n u_n-\n u|\d^{1-s_2} dx\bigg)^{\eta_1}\\ &\times &\bigg(\dyle\io|\n u_n-\n u|^{\g}\d^{\g(1-s_2)} dx\bigg)^{\frac{1-\eta_1}{\g}}, \\
        & \le & C \bigg(\dyle\io |\n u_n-\n u|\d^{1-s_2} dx \bigg)^{\eta_1}\to 0\mbox{  as   }n\to +\infty,
    \end{array}
    \right.
\end{equation}
where $\eta_1=\frac{\g-p\s}{p\s(\g-1)}$.\\
Since $p\s<\overline{m}_{s_1,s_2,q}$, then by repeating the same computations as in the first step, we can find $\b_1<2s_2-1$ such that
$$
\begin{array}{lll}
    \dyle\io|\n u_n-\n u|^{p}\d^{\b_1} dx &\le & C\bigg(\dyle\io|\n u_n-\n u|^{p\s}\d^{p\s(1-s_2)} dx\bigg)^{\frac{1}{\s}}\to 0\mbox{  as   }n\to +\infty.
\end{array}
$$
Hence , by using Proposition \ref{eq:Existence_with_weightel2}, we deduce that $v_n\to v$ strongly in $\mathbb{W}^{1,1}_0(\O)$ and then the continuity of $T$ follows.
\vspace{0.2cm}\\
$\bullet$ \textbf{\underline{Compactness} }\vspace{0.2cm}\\
-- Let $\{\varphi_n\}_n\subset H$ be such that
$\|\varphi_n\|_{\mathbb{W}^{1,1}_0(\O)}\le C \text{               and           } v_n=T(\varphi_n).$
Since $\{\varphi_n\}_n\subset H$, then 
$$\bigg(\io |\n \varphi_n|^{r}\, \d^{r(1-s_1)} dx\bigg)^{\frac{1}{r}}\le C\; \hbox{      for any       } n.$$

\noindent
Thus, as in the first step, there exists $\b<2s_1-1$ such that $\bigg\{|\n \varphi_n|^q \d^{\b}+\l f\bigg\}_n$ is bounded in $L^1(\O)$. Then,
according to Proposition \ref{eq:Existence_with_weightel2}, we deduce that, up to a subsequence, $u_n\to u$ strongly in $\mathbb{W}^{1,\a}_0(\O)$
for any $\a<\dfrac{N}{N-2s_1+1+\b}$.\\
\noindent
Recall that, from \eqref{eq:Strong_conv_uel2},  we have  proven that  $|\n u_n|\d^{1-s_2}\to| \n u|\d^{1-s_2}$ strongly in $L^{p\s}(\O)$ .\\
Thus by using the same arguments as in the first and second steps, we get the existence of  $\b_1<2s_2-1$ such that $(|\n u_n|^{p}+\mu g)\d^{\b_1}\to (|\n u|^{p}+\mu g)\d^{\b_1}$ strongly $L^{1}(\O)$.
Thanks to Proposition \ref{eq:Existence_with_weightel2}, we conclude that  $v_n\to v$ strongly in {$\mathbb{W}^{1,\a_1}_0(\O)$ for any $\a_1<\dfrac{N}{N-2s_2+1+\b_1}$}. In particular, $v_n\to v$ strongly in $\mathbb{W}^{1,1}_0(\O)$.  Hence, $T$ is compact.
\\
Therefore, by applying { Schauder's }fixed point Theorem, the
operator $T$ admits at least one fixed point $v\in H$ such that
$T(v)=v$. Hence, we get the existence of  $(u,v)\in
\mathbb{W}^{1,1}_0(\O)\times \mathbb{W}^{1,1}_0(\O)$ which is the
solution of System $(S)$ such that  $(u,v)\in
\mathbb{W}^{1,\theta_1}_0(\O,\d^{1-s_1})\times \mathbb{W}^{1,
\theta_2}_0(\O,\d^{1-s_1})$ for any
$\theta_1<\overline{m}_{s_1,s_2,q}$ and
$\theta_2<\widehat{\s}_{s_1,s_2,p}$.

Now, we briefly describe the main changes regarding  the set \eqref{eq:condi0012022}.\\
{\bf Second  case :} now, assume that the assumption \eqref{eq:condi0012022} is satisfied.\\
The main difference between the two sets is revealed when we show
that for every $\varphi\in H$, there exists $\beta <2s_1-1$ and
$\beta_1<2s_2-1$  such that $|\n \varphi|^q\d^{\b},\,{|\n
u|^p\d^{\b_1}\in L^1(\O)}$.

\noindent Since $s_1<\dfrac{q+1}{q+2}$ and
$m>\dfrac{1}{2s_1-q(1-s_1)}$, it follows that
$(q+1-s_1(q+2))m'<1$. Then, we can find $\b<2s_1-1$ such that
$(q(1-s_1)-\beta)m'<1$. Hence,  $\dyle \io \d^{(\beta-q(1-s_1))m'}
dx<\infty$ (respectively $\dyle\io \d^{(\beta_1-p(1-s_2))\s'}
dx<+\infty$). 
$\square$

\

\subsubsection{Case:  $p=q=1$}
Since the convexity argument cannot be applied in the case $p=q=1$, we will apply a different approach to establish the existence of a solution. 
{A natural step in this direction involves recasting the form System~$(S)$. Thus, by rescaling the space variable~$x$,  we obtain:}
\begin{equation}\label{eq:S1smalel2}
	\left\{
	\begin{array}{rclll}
		(-\Delta)^{s_1} u &= & \l\bigg(|\nabla v|+f\bigg) & \text{ in}&\Omega,\vspace{0.2cm} \\
		(-\Delta)^{s_2} v &= & \mu\bigg(|\nabla u|+g\bigg) & \text{ in}&\Omega, \\
		u=v &=& 0 &\hbox{  in}& \mathbb{R}^N\setminus\Omega.\\
	\end{array}%
	\right.
\end{equation}

\noindent
By using a suitable compactness result,  we get  the  following  existence result in this case.
\begin{thm}\label{eq:pq=1}  
 Assume that $(f,g)\in L^m(\O)^+\times L^ \s(\O)^+$ where $(m,\sigma)$ satisfies one of the following
    conditions:
    \begin{equation}\label{eq:FirstSet}
        \left\{
        \begin{array}{lll}
          2/3<s_1<1,\,\, 2/3 <s_2<1,\\
            1<m<\min\left\{\frac{N}{2s_1-1},\widehat{\s}_{s_2} \right\},\\ 1<\s<\min\left\{\frac{N}{2s_2-1},  \overline{m}_{s_1,s_2,1}\right\},
        \end{array}%
        \right.
    \end{equation}
    or
    \begin{equation}\label{eq:SecondSet}
        \left\{
        \begin{array}{lll}
      1/2<s_1<2/3,\,\, 1/2<s_2<2/3,\\ \frac{1}{3s_1-1}<m<\min\left\{\frac{N}{2s_1-1},{\widehat{\s}_{s_2}}  \right\},\\ \frac{1}{3s_2-1}<\s<\min\left\{\frac{N}{2s_2-1},  \overline{m}_{s_1,s_2,1}\right\},
        \end{array}%
        \right.
    \end{equation}
    or
    \begin{equation}\label{eq:ThirdSet}
        \left\{
        \begin{array}{lll}
        1/2<s_1<2/3,\,\, 2/3<s_2<1,\\ \frac{1}{3s_1-1}<m<\min\left\{\frac{N}{2s_1-1},{\widehat{\s}_{s_2}}   \right\},\\ 1<\s<\min\left\{\frac{N}{2s_2-1},  \overline{m}_{s_1,s_2,1}\right\},
        \end{array}%
        \right.
    \end{equation}
    or
    \begin{equation}\label{eq:FifthSet}
        m\ge \frac{N}{2s_1-1}\hspace{1mm}\,{\text{ and  }\;\; \frac{mN}{N+m(2s_2-1)}<\s<\frac{1}{s_2-s_1}},
    \end{equation}
    or {\begin{equation}\label{eq:SixthSet}
            \s \ge \frac{N}{2s_2-1} \; \text{ and  } \; m>\frac{\s N}{N+\s((2s_1-1)-N(s_2-s_1))}.\quad
    \end{equation}}
    Then,  there exist $\l^*>0$ and $\mu^*>0$ such that if $\l<\l^*$
    and $\mu<\mu^*$, System \eqref{eq:S1smalel2} has a solution
    $(u,v)$ such that $(u,v)\in
        \mathbb{W}^{1,\theta_1}_0(\O,\d^{\theta_1(1-s_1)})
        \times\mathbb{W}^{1,\theta_2}_0(\O,\d^{\theta_2(1-s_1)})$ for
    any $\theta_1<\overline{m}_{s_1,s_2,1}$ and
 $\theta_2<\widehat{\s}_{s_2}$.
\end{thm}
%
%

\bigbreak
\noindent
{\it Proof.}
    We give the proof under the assumption~\eqref{eq:FirstSet}. The others are similar.\\ By using  Schauder's fixed point Theorem, we get the existence of a solution $(u_n,v_n)$  to the following system
    \begin{equation}\label{eq:S1apr-}
        \left\{
        \begin{array}{rclll}
            (-\Delta)^{s_1} u_n &= & H_{1n} & \text{ in}&\Omega , \\
            (-\Delta)^{s_2} v_n &= & H_{2n} & \text{ in}&\Omega , \\
            u_n=v_n &=& 0 &\hbox{  in}& \mathbb{R}^N\setminus\Omega,\\
        \end{array}%
        \right.
    \end{equation}
    \noindent {where}  $$ H_{1n}:=\lambda\left(\dfrac{|\nabla v_n|}{1+\frac 1n |\n v_n|}+\dfrac{f}{1+\frac 1n f}\right)  \quad\text{ and }\quad
    H_{2n}:=\mu\left(\dfrac{|\nabla u_n|}{1+\frac 1n |\n u_n|}+\dfrac{g}{1+\frac 1n g}\right).
    $$

    \noindent

    By using the regularity result    in Theorem \ref{eq:regubeta1} with $a=1$, $s=s_1$ and $h=H_{1n}$, it follows that  
   \begin{equation}\label{eq:s00apr}
   	\left\||\n u_n|\d^{1-s_1}\right\|_{L^{\theta_1}(\O)} \leq  C\l \bigg(\left\||\n
   	v_n |\d^{1-s_1}\right\|_{L^{m}(\O)}+\|f\|_{L^{m}(\O)}\bigg),\qquad \forall \theta_1\le \dfrac{m N}{N-m(2s_1-1)}.
   \end{equation}
   {On the other hand, by using Theorem
   	\ref{eq:main_regularity_resultel2} with $a=1$ and $h=H_{2n}$, we
   	get
   	\begin{equation}\label{eq:s00aprmai}
   		\left\||\n v_n|\d^{1-s_1}\right\|_{L^{\theta_2}(\O)} \leq C\mu \bigg(\left\||\n u_n |\d^{1-s_1}\right\|_{L^{\s}(\O)}+\left\|g\right\|_{L^{\s}(\O)}\bigg)\qquad \forall\theta_2\le \dfrac{\s N}{N-\s(2s_2-1)}\\[9pt].
   \end{equation}}
   
   \noindent
   Then, we obtain
   \begin{equation}\label{eq:s01apr-}
   	\left\{
   	\begin{array}{lcc}
   		\left\||\n u_n|\d^{1-s_1}\right\|_{L^{\theta_1}(\O)}+ \left\||\n
   		v_n|\d^{1-s_1}\right\|_{L^{\theta_2}(\O)} \leq
   		\\
   		C\bigg(\mu \left\||\n u_n|\d^{1-s_1}\right\|_{L^{\s}(\O)}+\mu \|g\|_{L^{\s}(\O)} +\l \left\||\n
   		v_n|\d^{1-s_1}\right\|_{L^{m}(\O)}+\l \|f\|_{L^{m}(\O)}\bigg).
   	\end{array}
   	\right.
   \end{equation}
   Since {$\s<\overline{m}_{s_1,s_2,1}$} , we can
   choose $\theta_1$ close to
   {$\overline{m}_{s_1,s_2,1}$} such that
   $\s<\theta_1<\overline{m}_{s_1,s_2,1}$. Therefore, by using H\"older's
   inequality, we deduce that
   $$
   \left\||\n u_n |\d^{1-s_1}\right\|_{L^{\s}(\O)}\le C\left\||\n u_n
   |\d^{1-s_1}\right\|_{L^{\theta_1}(\O)}.
   $$
   In the same way, we  obtain
   $$
   \left\||\n v_n |\d^{1-s_1}\right\|_{L^{m}(\O)}\le  C\left\||\n v_n
   |\d^{1-s_1}\right\|_{L^{\theta_2}(\O)}.
   $$
    Hence, by choosing $\max\{\l,\mu\}<C_0$ such that $C_0$ is small enough, it follows that
{$$ \left\||\n
        u_n|\d^{1-s_1}\right\|_{L^{\theta_1}(\O)}+ \left\||\n
        v_n|\d^{1-s_1}\right\|_{L^{\theta_2}(\O)}\leq C(\O,C_0)\bigg(\l
        \|f\|_{L^{m}(\O)}+ \mu \|g\|_{L^{\s}(\O)}\bigg),
        $$}
    for every $\theta_1 <\overline{m}_{s_1,s_2,1}$ and $\theta_2 <{\widehat{\s}_{s_2}}$.\\
    As consequence, it follows that $\{(u_n,v_n)\}_n$ is bounded in
    $\mathbb{W}^{1,\theta_1}_0(\O,\d^{\theta_1(1-s_1)})\times
    \mathbb{W}^{1,\theta_2}_0(\O,\d^{\theta_2(1-s_1)})$ for all $\theta_1<\overline{m}_{s_1,s_2,1}$ and {$\theta_2<\widehat{\s}_{s_2}$}. \\
    Finally, the rest of the proof follows  by using the same compactness computations as in the proof of Theorem \ref{eq:existence_theorem_for_Sel2}.

\begin{remq}
The same existence {results hold} without any positivity
        assumptions on the data. $\square$
\end{remq}


        %

\subsection{Nonexistence results} 
In this subsection,  we will establish three  nonexistence results.
As we will see, those results highlight, in a certain sense,  the optimality of  the assumptions imposed on the data in Theorem \ref{eq:existence_theorem_for_Sel2}.  

\subsubsection{First nonexistence result}
\begin{thm}\label{eq:NonexisThm1}  
     Let  $(p,q)$ such that 
    {$1< q<\dfrac{s_1}{1-s_1}$ and $1< p<\dfrac{s_2}{1-s_1}$}. Assume that $(f,g) \in (L^{\infty}(\O)^+)^2$.  Then, there exist two positive constants $\Gamma(f)$ and $\Gamma(g)$ such that if
    $\l>\Gamma(f)$ or $\mu>\Gamma(g)$, System $(S)$ does not have
    any nonnegative solution.
\end{thm}
\bigbreak
\noindent
{\it Proof.}
    Let $\lambda, \mu >0$ be fixed such that System $(S)$ has a
    nonnegative solution $(u,v)$. Given $\theta\in
    \mathcal{C}_0^{\infty}(\Omega)$ with $\theta \gneqq 0$. Let $\phi_{\theta}\in \mathbb{W}_0^{s_1,2}(\Omega) \cap L^\infty(\O)$ and  $\varphi_{\theta}\in \mathbb{W}_0^{s_2,2}(\Omega) \cap L^\infty(\O)$ be the unique positive solutions to the following
    problems respectively
    \begin{equation} \label{eq:eqphi0-}
        \left\{
        \begin{array}{rclll}
            (-\D)^{s_1} \phi_{\theta} & = &  \theta & \quad \textup{ in } & \Omega\,, \\
            \phi_{\theta} & = &0 & \quad \textup{ in } & \mathbb{R}^N
            \setminus {\Omega},
        \end{array}
        \right.
    \end{equation}
    and
    \begin{equation} \label{eq:varphi0}
        \left\{
        \begin{array}{rclll}
            (-\D)^{s_2} \varphi_{\theta} & = &  \theta & \quad \textup{ in } & \Omega\,, \\
            \varphi_{\theta} & = &0 & \quad \textup{ in } & \mathbb{R}^N
            \setminus {\Omega}.
        \end{array}
        \right.
    \end{equation}
    Since $\theta$ is bounded, and according to \cite{RosOtonSerra}, we know
    that $\phi_\theta \backsimeq \d^{s_1}$ and $\varphi_\theta
    \backsimeq \d^{s_2}$. Moreover, let $\psi_{\theta}$ be  the unique solution
    to the following quasi-linear problem
    \begin{equation} \label{eq:psq-}
        \left\{
        \begin{aligned}
            -\text{div}( \varphi_{\theta} |\nabla \psi_{\theta}|^{\a-2} \nabla \psi_{\theta}) & = &\theta & \quad \textup{ in } &\Omega\,,\\
            \psi_{\theta} & = & 0 &\quad \textup{ on } &\partial \Omega,
        \end{aligned}
        \right.
    \end{equation}
    where $1<\a\le \min\{p,q\}$. {Let us mention} that the existence of $\psi_{\theta}$ follows by using an approximating argument.
   {Furthermore, it} is clear that $(\a-1)p'\le \a$ and $(\a-1)q'\le \a$.
    \\
    Now, let us use $\phi_\theta$ as a test function in the first equation of System $(S)$, we get
    \begin{equation}
        \left\{
        \begin{array}{lll}
            \dyle\int_\Omega \left(|\n v|^q+\l f\right)\phi_\theta\, dx &= & \dyle\int_\Omega  \theta u\, dx,\\
            & \le & \dyle\int_\Omega  \varphi_{\theta} |\nabla \psi_{\theta}|^{\a-1} |\nabla u|\, dx.
        \end{array}
        \right.
    \end{equation}
    Therefore,  {by }using H\"older's and Young's inequalities, it follows that
    \begin{equation}\label{eq:save1}
        \left\{
        \begin{array}{lll}
            \dyle\int_\Omega \left(|\n v|^q+\l f\right)\phi_\theta\, dx &\le &
            \bigg(\dyle\io |\nabla u|^p \varphi_\theta\, dx \bigg)^{\frac 1p}\bigg(\dyle\io |\nabla \psi_\theta|^{p'(\a-1)} \varphi_\theta\, dx \bigg)^{\frac{1}{p'}},\\
            &\leq &\varepsilon \dyle\io |\nabla u|^p \varphi_\theta\, dx+ {C_1}'(\varepsilon)\dyle\io |\nabla \psi_\theta|^{p'(\a-1)} \varphi_\theta\, dx. \\
        \end{array}
        \right.
    \end{equation}
    \noindent
    Analogously, using $\varphi_\theta$  as a test function in the
    second equation of System $(S)$, we get
    \begin{equation}\label{eq:save100}
        \begin{array}{lll}
            \dyle\int_\Omega \left(|\n u|^p +\mu g\right)\varphi_\theta\, dx &\le &\varepsilon \dyle\io |\nabla v|^{q} \varphi_\theta\, dx+
            C'_2(\varepsilon)\io |\nabla \psi_\theta|^{q'(\a-1)}\varphi_\theta\, dx.
        \end{array}
    \end{equation}

    \noindent
    Combining \eqref{eq:save1} and \eqref{eq:save100}, it holds that
    \begin{equation}\label{eq:s001}
        \left\{
        \begin{array}{lll}
            \dyle \int_\Omega \left(|\n v|^q+\l f\right)\phi_\theta\, dx &\leq& \varepsilon^2\dyle\io |\nabla v|^{q} \varphi_\theta\, dx+{\e C'_2(\e)}\io |\nabla \psi_\theta|^{q'(\a-1)}\varphi_\theta\, dx\\ &+&{C_1}'(\e)\dyle\io |\nabla \psi_\theta|^{p'(\a-1)}\varphi_\theta\, dx.
        \end{array}
        \right.
    \end{equation}

    \noindent
    By using the fact that,  $\varphi_\theta\backsimeq \d^{s_2}$ and
    $\phi_\theta\backsimeq \d^{s_1}$ in the previous inequality, we obtain
    \begin{equation}\label{eq:s002}
        \left\{
        \begin{array}{lll}
            \dyle \int_\Omega \left(|\n v|^q+\l f\right)\d^{s_1}\, dx &\leq \varepsilon^2 C\dyle\io |\nabla v|^{q} \d^{s_2}\, dx+\overline{C_1}(\e)\io |\nabla \psi_\theta|^{q'(\a-1)}\d^{s_2}\, dx\\ &+ \overline{C_2} (\e)\dyle\io |\nabla \psi_\theta|^{p'(\a-1)}\d^{s_2}\, dx.
        \end{array}
        \right.
    \end{equation}
    Since $\O$ is bounded  and $s_1<s_2$, then $\d^{s_2}(x)\leq C(\O)\d^{s_1}(x)$ for any $x\in \O$.\\
    Thus, from the inequality \eqref{eq:s002}, it follows that
    \begin{equation}\label{eq:s0022}
        \left\{
        \begin{array}{lll}
            \dyle \int_\Omega \left(|\n v|^q+\l f\right)\d^{s_1}\, dx &\leq \varepsilon^2 C\dyle\io |\nabla v|^{q} \d^{s_1}\, dx+\overline{C_1}(\e)\io |\nabla \psi_\theta|^{q'(\a-1)}\d^{s_2}\, dx\\ &+ \overline{C_2}(\e)\dyle\io |\nabla \psi_\theta|^{p'(\a-1)}\d^{s_2}\, dx.
        \end{array}
        \right.
    \end{equation}
    Hence, by choosing $0<\varepsilon <<\sqrt{\frac{1}{C}}$, we get
    \begin{equation}
        \l\int_\Omega  f \d^{s_1}\, dx \leq \dyle  \overline{C_1}(\e)\io |\nabla
        \psi_\theta|^{p'(\a-1)} \d^{s_2}\, dx +\overline{C_2}(\e)\io |\nabla
        \psi_\theta|^{q'(\a-1)} \d^{s_2}\, dx.
    \end{equation}

    \noindent       Otherwise, we obtain
    \begin{equation}\label{eq:lastineq}
        \left\{
        \begin{array}{lll}
            \l \dyle\int_\Omega  f \phi_\theta\, dx &\leq& C_3(\e)\dyle  \io |\nabla
            \psi_\theta|^{p'(\a-1)} \varphi_\theta\, dx +C_4(\e)\io |\nabla
            \psi_\theta|^{q'(\a-1)} \varphi_\theta\, dx,\\
            &\leq& C(\e)\left(\dyle  \io |\nabla \psi_\theta|^{p'(\a-1)} \varphi_\theta\, dx +\io |\nabla
            \psi_\theta|^{q'(\a-1)} \varphi_\theta\, dx\right),
        \end{array}
        \right.
    \end{equation}
    where $C(\e)=\max\big\{C_3(\e),C_4(\e)\big\}$.\\

    \noindent
    Now, let us set
    \begin{equation}\label{eq:equat1}
        F(\theta):=\dyle \io |\nabla \psi_\theta|^{p'(\a-1)} \varphi_\theta\, dx +\io |\nabla \psi_\theta|^{q'(\a-1)} \varphi_\theta\, dx.
    \end{equation}
    Then,  $\l\le \Gamma(f)=C(\varepsilon) \l^*(f)$ where
    \begin{equation*}\label{eq:alpha-etoile-}
        \l ^*(f):=\inf_{\{0\lneqq\theta\in
            \mathcal{C}^{\infty}_0(\Omega)\}}\left\{F(\theta)
        \;;\;\int_{\Omega}f\phi_\theta\, dx=1\right\}.
    \end{equation*}
    \noindent Conversely, if $\l >\Gamma(f)$, System $(S)$ does not have any nonnegative solution. By a similar way, we get the result for $ \mu$.

\subsubsection{Second nonexistence result}

\begin{thm}\label{eq:NonExisThm2-} 
    The following two assertions hold, even for small  $\lambda$ and $\mu$, and regardless of the regularity of  $f$ and $g$.\\
    \noindent\text{\em 1) } Assume that $\max\{p,q\}\ge \dfrac{1+s_2}{1-s_1}$. 
    Then, System~$(S)$ does not have any nonnegative solution
    $(u,v)\in \mathbb{W}^{1,1}_0(\Omega)\times \mathbb{W}^{1,1}_0(\O)$ such that
    $$\dint_\Omega |\n u|^p\d^{s_2}\, dx<+\infty\quad \text{ and  }\quad \dint_\Omega |\n v|^q\d^{s_1}\, dx <+\infty. $$
    \noindent\text{\em 2) } Assume that $\max\{p,q\}\ge \dfrac{1}{1-s_1}$. Then,
    System~$(S)$ does not have any nonnegative solution $(u,v)\in \mathbb{W}^{1,p}_0(\Omega)\times
    \mathbb{W}^{1,q}_0(\O)$.
\end{thm}

\bigbreak
\noindent
{\it Proof.}
    \
    \

    \

    \noindent {\bf{First case:}} Let us start by proving the first case. Without loss of generality, we can suppose that $p\ge q$ and  $p\ge \dfrac{1+s_2}{1-s_1}$.
    \noindent Assume by contradiction that System $(S)$ has a solution $(u,v)$ such that $(u,v)\in
    \mathbb{W}^{1,1}_0(\Omega)\times \mathbb{W}^{1,1}_0(\O)$ and
    $$\dint_\Omega |\n u|^p\d^{s_2}\, dx<+\infty\quad \text{ and  }\;\; \dint_\Omega |\n v|^q\d^{s_1}\, dx <+\infty. $$

    \smallbreak

    \noindent Since $s_2<p-1$, by  using  the weighted Hardy's
    inequality (Theorem \ref{eq:HardyInequalityel2} {with $\rho=s_2$}), we obtain
    \begin{equation}
        \dint_\Omega \frac{u^p}{\d^p}\d^{s_2}\,  dx\le \dint_\Omega |\nabla
        u|^p \d^{s_2}\,  dx<+\infty.\end{equation} From the maximum principle proved in  \cite{BM}, we have $u\ge C
    \d^{s_1}$. Hence,
    $$\int_\Omega \dfrac{1}{\d^{p(1-s_1)-s_2}} dx<+\infty.$$
    Thus,  $p<\dfrac{1+s_2}{1-s_1}$ which is contradiction.\\
    \noindent {\bf{Second case:} } Now, we prove  the second case. Suppose that
    System $(S)$ has a nonnegative solution $(u,v)$ with $(u,v)\in
    \mathbb{W}^{1,p}_0(\Omega)\times \mathbb{W}^{1,q}_0(\O)$. Then, by the same reasoning as in first case, by using Hardy's  inequality (Theorem \ref{eq:HardyInequalityel2} with $\rho=0$), it follows that
    $$
    \dint_\Omega\frac{u^p}{\d^p}\,  dx\le \dint_\Omega |\nabla u|^p\,  dx<+\infty.
    $$
    Furthermore, by using the fact that $u\ge C \d^{s_1}$, it holds that
    $$\int_\Omega \dfrac{1}{\d^{p(1-s_1)}} dx<+\infty.$$
    Hence, $p<\dfrac{1}{1-s_1}$ which is a contradiction. This completes the proof.

\subsubsection{Third nonexistence result}
In this paragraph, we state an optimality result related to the assumptions on $f$ and $g$.

\begin{thm}\label{eq:Nonexistence3}   
Assume that $p,q\geq1$ 
Then, there exists $(m,\s)$
	satisfying
    \begin{equation}\label{eq:noncondim-}
        p\s>\dfrac{N}{N-m(2s_1-1)}\;\mbox{  and  }\;qm>\dfrac{ N}{N-\s(2s_2-1)},
    \end{equation}
     and there exists $(f,g)\in L^m(\O)\times L^\s(\O)$ with
    $f, g\gneq 0$, such that  System $(S)$ does not have any positive solution, regardless of how small $\l$ and $\mu$ are.   
    \end{thm} \noindent 
 \textbf{\noindent  Comments.}
	Before giving the proof of Theorem \ref{eq:Nonexistence3}, let us illustrate Assumption \eqref{eq:noncondim-}.  More specifically, for $s_1=\frac34$, $s_2=\frac{9}{10}$, $N=3$ and fixed $(p,q)$, let us exhibit some  explicit examples of the values of $(m,\sigma)$ where the assumption is fulfilled.  
	\begin{itemize}
   \item[(i)] $p=q=1$ and one of the following  $\begin{cases}
  \sigma \in[1,\frac{25}{8}], \;\;  m \in(6,+\infty)
   	\\
    \sigma \in\left(\frac{15}{4},+\infty\right), \;\; m \in\left[1, \frac{6(\sigma-1)}{\sigma}\right)\cup   (6,+\infty) 
   	\\
\sigma \in\left(\frac{25}{8}, \frac{15}{4}\right),\;\; m \in\left(\frac{15}{15-4 \sigma},+\infty\right)
   	\\
	\sigma \in\left(\frac{27}{20}, \frac{11}{4}\right], \;\;  m \in\left(\frac{15}{15-4 \sigma}, \frac{6(\sigma-1)}{\sigma}\right).\end{cases}  $
   	\item[(ii)]  $p=1$, $q=2$ and one of the following $
   	\begin{cases}
   		\sigma \in[1,\frac{55}{16}]\cup  \left(\frac{15}{4},+\infty\right), \;\;  m \in(6,+\infty)
   		\\
   		\sigma \in\left(\frac{6}{5}, \frac{15}{8}\right) \cup\left(\frac{15}{4},+\infty\right),\;\;  m \in\left[1, \frac{6(\sigma-1)}{\sigma}\right)\\
   		\sigma=\frac{15}{8},\;\; m\in (1,\frac{14}{5})\\
   	\sigma \in\left(\frac{55}{16}, \frac{15}{4}\right),\;\; m \in\left(\frac{15}{2(15-4 \sigma)},+\infty\right)
   		\\
   	\sigma \in\left(\frac{15}{8}, \frac{33}{10}\right),\;\; m \in\left(\frac{15}{2(15-4 \sigma)},\frac{6(\sigma-1)}{\sigma}\right).
   	\end{cases}
   $
  	\item[(iii)]  $p=2$, $q=1$ and one of the subsequent
   $\left\{\begin{array}{l}
   \sigma=1,\;\;  m \in\left(\frac{15}{11}, 3\right)
   	\\
\sigma \in [1,\frac{25}{8}]\cup \left(\frac{15}{4},+\infty\right),\;\;  m \in\left(6,+\infty\right) \\   	
   		 \sigma \in\left(\frac{25}{8},\frac{15}{4}\right),\;\;  m \in\left[\frac{15}{15-4 \sigma},+\infty\right)
   	\\
 \sigma \in\left(\frac{15}{4},+\infty\right),\;\;  m \in\left[1, \frac{3(2 \sigma-1)}{\sigma}\right) \\
 \sigma \in(1,3),\;\;  m \in\left(\frac{15}{15-4 \sigma}, \frac{3(2 \sigma-1)}{\sigma}\right).
   \end{array}\right.$
   \item[(iv)]  
  $ p=q=2$ and one of the subsequent $\left\{\begin{array}{l}
    \sigma=1, \;\;  m \in [1,3) 		
   	\\
   	\sigma \in [1,\frac{55}{16}]\cup \left(\frac{15}{4},+\infty\right),\;\;  m \in\left(6,+\infty\right)
   	\\
   	\sigma \in (1,\frac{15}{8})\cup \left(\frac{15}{4},+\infty\right),\;\;  m \in\left[1,\frac{3(2\sigma-1)}{\sigma}\right)
   	\\
   	 \sigma=\frac{15}{8}, \;\; m \in (1,\frac{22}{5}) 		
   	\\
    \sigma \in\left(\frac{15}{8}, \frac{27}{8}\right] , \;\;  m \in\left(\frac{3(2 \sigma-1)}{\sigma},\frac{15}{2(15-4 \sigma)}\right)
    	\\
    \sigma \in\left(\frac{55}{16}, \frac{15}{4}\right) , \;\;  m \in\left(\frac{15}{2(15-4 \sigma)},+\infty\right).
   \end{array}\right.$
    \item[(v)]  $\forall p, q \geq 1$,  \eqref{eq:noncondim-} is satisfied  $ \;\forall \sigma \in[1,3] \cup\left(\frac{15}{4},+\infty\right),\;\; \forall m \in(6,+\infty)$.    
     \item[(vi)] $\forall p \geq 1$, if $ q=1$  then \eqref{eq:noncondim-} holds   $\;\forall \sigma \in\left(\frac{25}{8}, \frac{15}{4}\right),\;\; \forall m \in\left(\frac{15}{15-4 \sigma},+\infty\right) $.
      \item[(vii)] $\forall p \geq 1,$ if  $ q=2$  then \eqref{eq:noncondim-} is fulfilled   $\;\forall \sigma \in\left(\frac{55}{16}, \frac{15}{4}\right),\;\; \forall m \in\left(\frac{15}{2(15-4 \sigma)}, +\infty\right)$.
	\end{itemize}

\noindent{\bf Proof of Theorem \ref{eq:Nonexistence3}.} 
    {Let us construct $f$ and $g$ such that, under the assumptions of Theorem \ref{eq:NonexisThm1}, we get $\Gamma(f)=0$ or $\Gamma(g)=0$.\\
        Without loss of generalities, we suppose that {$\O=B_1(0)$} and $m,\s<N$. Let $\e>0$, that will be chosen later,
        and define $f(x)=\frac{1}{|x|^{\frac{N-\e}{m}}}$. Then,
        $f\in L^m(\O)$. Assume by contradiction that System $(S)$ has a solution $(u,v)\in \mathbb{W}^{1,1}_0(\O)\times \mathbb{W}^{1,1}_0(\O)$ such that $(u,v)\in
        \mathbb{W}^{1,p}_{loc}(\O)\times \mathbb{W}^{1,q}_{loc}(\O)$.
        Then, by following the calculations done to prove Theorem \ref{eq:NonexisThm1}, we
        deduce that
\begin{equation}\label{eq:nancy11}
            \l\int_\Omega  f\phi_\theta\, dx \leq \dyle C\io (|\nabla
            \psi_\theta|^{p'(\a-1)}+ |\nabla
            \psi_\theta|^{q'(\a-1)} ) \varphi_\theta\, dx.
    \end{equation}
        Moreover, let us consider the {following} function 
        $$ \theta(x)= \left\{
        \begin{array}{lll}
            \dfrac{1}{|x|^a}& \mbox{  if  }& |x|\le \frac 14\\
            (1-|x|)^\g & \mbox{  if  }& \frac 12\le |x|\le 1,
        \end{array}
        \right.
        $$
        where {$\gamma >1$} and {$2s_1<a<N+2s_1$} satisfying $$ \left\{
        \begin{array}{lll}
            N+2s_1-\frac{N}{m}<a<& N+2s_1-p'(2s_2-1), \\
            N+2s_1-\frac{N}{\s}<a< &N+2s_1-q'(2s_2-1).
        \end{array}
        \right.
        $$
        {Let us mention} that the existence of $a$ follows using the assumption
        \eqref{eq:noncondim-}.
        \\
{Thus, by denoting $G(x):=(|\nabla
        \psi_\theta(x)|^{p'(\a-1)}+ |\nabla
        \psi_\theta(x)|^{q'(\a-1)} ) \varphi_\theta(x)$, we have
\begin{equation}
                \io G(x) dx =\int_{B_{\frac14}(0)}G(x) dx+\int_{\{\frac14<|x|<\frac12\}}G(x) dx+\int_{\{\frac12<|x|<1\}}G(x) dx.
        \end{equation}
}
        \noindent By taking into consideration that $\phi_{\theta}, \varphi_{\theta}$ solve Problems \eqref{eq:eqphi0-} and \eqref{eq:varphi0} respectively,
we obtain that $\phi_{\theta}\simeq
            |x|^{2s_1-a}$ and $\varphi_{\theta}\simeq
                |x|^{2s_2-a}$ in a small ball $B_{r_0}(0)\subset\subset B_{1}(0)$.
        \smallbreak
        \noindent Regarding the function $\psi_{\theta}$ which is the solution to
        Problem \eqref{eq:psq-}, choosing $1<\a<2s_2$, we can prove that in a
        small ball $B_{r_1}(0)\subset\subset B_{1}(0)$, we have
        \begin{equation}\label{eq:ponct1}
            \psi_{\theta}(x)\simeq |x|^{-\frac{2s_2-\a}{\a-1}}\; \mbox{  and }\;|\n
            \psi_\theta|\le C|x|^{-\frac{2s_2-1}{\a-1}}.
        \end{equation}
      Furthermore, by taking into consideration the estimates
    \eqref{eq:nancy11} and \eqref{eq:ponct1} and by choosing
    $\e$ small enough such that $0<\e<<N-m(N-a+2s_1)$, we deduce
    that
 $$
    \dyle\io f\phi_\theta  dx\ge C\dyle\int_{B_{\frac14}(0)}
    \frac{1}{|x|^{\frac{N-\e}{m}+a-2s_1}} dx=+\infty.
    $$
On the other hand, it is clear that, for $\frac14<|x|<\frac12$,
        \begin{equation}
            \int_{\{\frac14<|x|<\frac12\}} (|\nabla
            \psi_\theta|^{p'(\a-1)}+ |\nabla
            \psi_\theta|^{q'(\a-1)} ) \varphi_\theta\, dx<C.
    \end{equation}
Now, for any $x\in B_{\frac14}(0)$, we have
\begin{equation}
            \int_{B_{\frac14}(0)} (|\nabla
            \psi_\theta|^{p'(\a-1)}+ |\nabla
            \psi_\theta|^{q'(\a-1)} ) \varphi_\theta\, dx\le   C \dyle \io\left(
            \frac{1}{|x|^{(2s_2-1)p'+a-2s_2}}+\frac{1}{|x|^{(2s_2-1)q'+a-2s_2}}\right ) dx<+\iy.
    \end{equation}
With respect to the integral over the set $\{\frac12<|x|<1\}$, in a similar way,
we find that
        $$\int_{\{\frac12<|x|<1\}} (|\nabla
            \psi_\theta|^{p'(\a-1)}+ |\nabla
            \psi_\theta|^{q'(\a-1)} ) \varphi_\theta\, dx<C.
$$

    Therefore, using a suitable approximation argument, we conclude that
$\Gamma(f)=0$. Hence, the nonexistence result follows.} \hfill$\square$

%
%

%

\medbreak
To conclude this paper, it is worth mentioning that by using the same approach as in 
Section \ref{eq:Existenc_Nonexistence_S}, we can establish  existence results for the following systems.
\smallbreak
\noindent$\bullet$ {\bf Systems with drift:}
\begin{equation} 
	\label{eq:SDR--} 
    \left\{
    \begin{array}{rclll}
        (-\Delta)^{s_1} u+B_1(x).\n u &= & |\nabla v|^{q}+\l f & \text{ in}&\Omega , \\
        (-\Delta)^{s_2} v+B_2(x).\n v &= & |\nabla u|^{p}+\mu g & \text{ in}&\Omega , \\
        u=v &=& 0 &\hbox{  in}& \mathbb{R}^N\setminus\Omega,\\
    \end{array}%
    \right.
\end{equation}
where $B_i$ is a vector fields satisfying $|B_i|\in
L^{\gamma_i}(\O)$ with  $\gamma_i>\frac{N}{2s_i-1}$.
\smallbreak
\noindent$\bullet$ {\bf Systems  where at least one of the two gradients acts as an absorption term:}
        \begin{equation}\label{eq:PSel2}
            \left\{
            \begin{array}{rclll}
                (-\Delta)^{s_1} u +|\nabla v|^{q}&=& \l f & \textnormal{ in } &\Omega , \\
                (-\Delta)^{s_2} v \pm|\nabla u|^{p}&=& \mu g & \textnormal{ in }&\Omega , \\
                u=v &=& 0 &\textnormal{  in }& \mathbb{R}^N\setminus\Omega.
            \end{array}%
            \right.
        \end{equation}%
        \hfill$\blacksquare$

\end{document}